\documentclass[12pt]{amsart}
\usepackage{amssymb,amsfonts,amsthm,amsmath,amscd}

\title[Isocohomological Property and Dehn Functions]{The Isocohomological Property, Higher Dehn Functions, and Relatively Hyperbolic Groups}

\author{Ronghui Ji}
\address{Department of Mathematical Sciences\\
Indiana University-Purdue University, Indianapolis\\ 
402 N. Blackford Street\\
Indianapolis, IN 46202}
\email{ronji@math.iupui.edu}

\author{Bobby Ramsey}
\email{bramsey@math.iupui.edu}

\subjclass[2000]{Primary 20F65; Secondary 58B34 }
\keywords{polynomial cohomology, higher Dehn functions, relatively hyperbolic groups, Novikov conjecture}

\date{}

\newtheorem{thm}{Theorem}[section]
\newtheorem{defn}[thm]{Definition}

\newtheorem{cor}[thm]{Corollary}
\newtheorem{lem}[thm]{Lemma}

\newcommand{\Z}{{\mathbb{Z}}}
\newcommand{\R}{{\mathbb{R}}}
\newcommand{\C}{{\mathbb{C}}}

\newcommand{\N}{{\mathbb{N}}}

\newcommand{\SnX}[1][]{{\mathcal{S}_{#1} X}}
\newcommand{\CnX}[1][]{{C_{#1} X}}
\newcommand{\SG}{\mathcal{S}G}
\newcommand{\CG}{{\mathbb{C}G}}

\newcommand{\btensor}{{\hat\otimes}}

\newcommand{\isom}{{\,\cong\,}}
\newcommand{\bHom}{\operatorname{bHom}}
\newcommand{\Hom}{\operatorname{Hom}}
\newcommand{\Ext}{\operatorname{Ext}}

\allowdisplaybreaks[4]

\begin{document}

\begin{abstract}
The property that the polynomial cohomology with coefficients of a finitely generated discrete group is canonically isomorphic to the group cohomology is called the (weak) isocohomological property for the group. In the case when a group is of type $HF^\infty$, i.e. that has a classifying space with the homotopy type of a polyhedral complex with finitely many cells in each dimension, we show that the isocohomological property is geometric and is equivalent to the property that the universal cover of 
the classifying space has polynomially bounded
higher Dehn functions.  If a group is hyperbolic relative to a collection of subgroups, each of which is polynomially combable, respectively $HF^\infty$ and isocohomological, then we show that the group itself has these respective properties.  Combining with the results of Connes-Moscovici and Dru{\c{t}}u-Sapir we conclude that a group satisfies the strong Novikov conjecture if it is hyperbolic relative to subgroups which are of property RD, of type $HF^\infty$ and
isocohomological.
\end{abstract}   

\maketitle


\section{Introduction}
Given a finitely presented group $G$ with classifying space $BG$, which is a $K(G, 1)$ space, and a compact oriented smooth manifold $M$, together with a continuous map $\varphi: M \longrightarrow BG$, the higher signatures of the pair $(M, \varphi)$ are defined to be $(L(M)\cdot \varphi^*(\xi), [M])$, where $\xi$ is any class in $H^*(BG, Q)$ and $L(M)$ denotes the total Hirzebruch $L$-class of $M$. The Novikov conjecture \cite{Novikov} states that the higher signatures defined above are homotopy invariants of the pair $(M, \varphi)$. The validity of this conjecture has been established, for many classes of groups by variety of techniques (see \cite{HK, Rosenberg, yu2}). Using Kasparov's $KK$-theory \cite{Kas1}, the homotopy invariance of the higher signatures is a consequence of the rational injectivity of the K-theoretical assembly map from $K$-homology of the classifying space $BG$ to the $K$-theory of the reduced group $C^*$-algebra $C_r^*(G)$. The conjecture that this assembly map is always rationally injective is called the strong Novikov conjecture \cite{BCH}.  

The early work of Connes and Moscovici \cite{CM} shows that if a finitely generated discrete group $G$ satisfies the following two conditions then $G$ satisfies the strong Novikov conjecture.  

(1) G has the Rapid-Decay property of Jolissaint \cite{jol}.

(2) The polynomial cohomology $HP^*(G; \C)$ is surjective to $H^*(G; \C)$ for the homomorphism induced by the inclusion.

Here $HP^*(G; \C)$ is the polynomial growth cohomology of the group $G$. This cohomology theory was proposed by Connes and Moscovici in \cite{CM} and formalized in \cite{Ji}. Specifically, let $|g|$ denote the word
length of an element $g$ in the group $G$, with respect to some fixed finite generating set. An $n$-cochain $c$ on $G$ is of
polynomial growth if for all elements $g_0$, $g_1$, ..., $g_n$ in $G$, 
$|c(g_0, g_1, ..., g_n)| \le P(|g_0| + |g_1| + \cdots + |g_n|)$ for some polynomial $P$ depending only on $c$.  The collection of polynomial
growth cochains forms a sub-cochain complex of the usual group cochain complex with complex coefficients.  The cohomology of this subcomplex is
called the polynomial cohomology of $G$ with coefficients $\C$.  The inclusion of the subcomplex into the full complex
induces the comparison homomorphism $HP^*(G; \C) \rightarrow H^*(G; \C)$.
More generally, one can define polynomial cohomology with coefficients in the category of Fr\'echet spaces \cite{Ji} as follows. Let $S_1(G) = \{f: G \longrightarrow \C | \sum_{g\in G}|f(g)|(1 + |g|)^k < \infty, k >0\}$. $S_1(G)$ is clearly a topological algebra \cite{jol2} 
with the natural Fr\'echet topology given by the seminorms $\|f\|_k = \sum_{g\in G}|f(g)|(1 + |g|)^k $. Let $V$ be any 
Fr\'echet space upon which $S_1(G)$ acts continuously. The polynomial cohomology $HP^*(G; V)$ of $G$ with coefficients in $V$ is defined to be $\Ext_{S_1(G)}^*(\C, V)$ in the category of continuous $S_1(G)$-modules. Note that in this category, a projective resolution of $\C$ over $S_1(G)$ must be endowed with a continuous $\C$-linear splitting. 
Since word hyperbolic groups satisfy both conditions \cite{Gr, Ha}, Connes and Moscovici conclude that word hyperbolic groups satisfy the Novikov conjecture. Until recently \cite{Ji-Ogle, meyer, Ogle, Ramsey} there has been little progress in verifying the Novikov conjecture using this method.
The main difficulty is the verification of condition (2). For condition (1) there has been much progress made in a variety of cases \cite{Ch1, CR, DS, La, RRS}.

A finitely generated group $G$ has the {\it{weak isocohomological property}} if for every coefficient module $V$,
the comparison homomorphism $HP^*(G; V) \rightarrow H^*(G;V)$ is an isomorphism.  
The term `isocohomological' is taken from Meyer \cite{meyer2}, where it describes a homomorphism between two bornological
algebras.  What is meant here by `$G$ has the weak isocohomological property', is a weakened version of 
what Meyer refers to as `the embedding $\C[G] \rightarrow S_1(G)$ is isocohomological'.  Hereafter we drop the 
`weak' adjective, and refer to this as the isocohomological property of the finitely generated group $G$.

The first breakthrough in this subject was by Meyer \cite{meyer} and Ogle \cite{Ogle} who independently proved that any polynomially combable group has the isocohomological property.  Not only so, Ogle \cite{Ogle} proves that if $G$ is of type $HF^\infty$, (i.e. that $G$ has a classifying space with the homotopy type of a polyhedral complex with finitely many cells in each dimension) and satisfies polynomial growth for his version of higher Dehn functions, then $G$ is
isocohomological. Ogle's higher Dehn functions are defined to be any contracting homotopy of the topological projective resolution of $\C$ over $S_1(G)$ obtained from the universal cover of a classifying space that has finite many cells in each dimension.  As such the relationship between Ogle's higher Dehn functions and the usual higher Dehn functions, as first studied in \cite{alonso}, is not clear.
We remark that a polyhedral complex is an analogue of a simplicial complex, without the rigidity
that for each dimension $n$, every $n$-cell has a fixed number of faces \cite{alek}.

In this paper we introduce the concept of weighted fillings for $n$-boundaries over $G$ and define the higher weighted Dehn functions for an $HF^\infty$ group $G$ by using these weighted fillings. 
Our version is equivalent to Ogle's in the sense that when our weighted Dehn functions are of polynomial growth, so are Ogle's complex admits polynomially bounded higher Dehn functions, and vice-versa.  We also show that the weighted higher Dehn functions having polynomial growth is equivalent to the usual higher Dehn functions, as studied by Gersten \cite{ger}, having polynomial growth. We further prove that for $HF^\infty$ groups, the isocohomological property 
is equivalent to the usual higher Dehn functions of the group having polynomial growth. One of the key ingredients in the proof is a technique used by Mineyev \cite{min}, in which he proved that a group is word hyperbolic if and only if the degree two bounded cohomology with coefficients for the group is surjective onto the usual group cohomology under the map induced by the inclusion. Since the Dehn functions are equivalent for quasiisometric groups, one asks whether the isocohomological property is preserved for quasiisometric groups. We show that this is indeed the case among $HF^\infty$ groups. As a consequence all groups of polynomial growth have polynomially bounded higher Dehn functions. For a finitely generated nilpotent group it is known that the first Dehn function is polynomially bounded \cite{ger2}. The details are given in section 2.  We note that another class of $HF^\infty$ groups that have all higher Dehn functions of polynomial growth is the class of groups that possess a polynomial combing \cite{ger}. 

Recently Osin \cite{Os} and Dadarlat-Guentner \cite{DG}, respectively, proved that if a group is relatively hyperbolic with respect to subgroups that have finite asymptotic dimensions, respectively are coarsely embeddable in Hilbert spaces, then the group itself has finite asymptotic dimension, respectively is coarsely embeddable in a Hilbert space. These two results imply the validity of the Novikov conjecture for both classes of groups via the method of coarse geometry \cite{yu1,yu2}. One asks whether or not Connes-Moscovici's method also applies to relatively hyperbolic groups. For this purpose we first construct a polynomial combing for groups that are relatively hyperbolic to subgroups which are polynomially combable.  Thus, after observing the work of Dru{\c{t}}u-Sapir \cite{DS} which states that a group relatively hyperbolic to subgroups of property RD itself has property RD, Connes-Moscovici's method does work for this class of relatively hyperbolic groups.
We remark that by `relatively hyperbolic' we mean relatively hyperbolic with the bounded coset penetration property \cite{farb}.  The above quoted theorem of Dru{\c{t}}u-Sapir was discovered by Chatterji-Ruane in the case that 
the subgroups are of polynomial growth \cite{CR}.


Since polynomially combable groups are $HF^\infty$ groups and have all higher Dehn functions of polynomial growth, it is desirable to extend the results in section 3 to the larger class of groups which contains all $HF^\infty$ groups of isocohomological property. By using the method of the combing constructions in section 3, we prove that if a group is relatively hyperbolic to subgroups that are of type $HF^\infty$ and are isocohomological, then the group itself is of type $HF^\infty$ and is isocohomological. This is done by constructing a classifying space for the group of $HF^\infty$-type and by estimating the growth of the higher Dehn functions of the group in terms of those of the subgroups. This will be done in the final section. Again by the result of Dru{\c{t}}u-Sapir \cite{DS}, we conclude that when the subgroups are of property RD, of type $HF^\infty$ and isocohomological, then the group satisfies the strong Novikov conjecture. 

Finally, we would like to thank Crighton Ogle, Mark Sapir and Jonathan Rosenberg for some 
helpful comments and discussions.

\section{Characterization of polynomial growth cohomology}
Recall from the introduction that the polynomial cohomology of a finitely generated discrete group $G$ with coefficients in a Fr{\'e}chet space $V$ is defined as $\Ext_{S_1(G)}^*(\C, V)$. To calculate such cohomology groups one needs a topologically projective resolution with continuous $\C$-linear splittings \cite{Taylor} of $\C$ over $S_1(G)$: 
\[ 0 \leftarrow \C \leftarrow P_0 \leftarrow P_1 \leftarrow \cdots \leftarrow P_n \leftarrow \cdots\ \ \ (*)\]
 Then $\Ext_{S_1(G)}^*(\C, V)$ is the homology group of the complex of continuous module homomorphisms: 
 \[ 0 \rightarrow \Hom_{S_1(G)}(P_0, V) \rightarrow \Hom_{S_1(G)}(P_1, V)\rightarrow \cdots \rightarrow \Hom_{S_1(G)}(P_n, V) \rightarrow \cdots \ \ \ (**)\]
The usual topological bar resolution \cite{Ji} is such a resolution but it is infinite dimensional after taking $\Hom$ in (**). In this section we wish to find a topological resolution (*) so that (**) becomes finite dimensional for each degree under certain assumptions for the group.

Let $X$ be a polyhedral complex, with $n$-skeleton denoted by $X^{(n)}$.  A collection of $(n+1)$-cells $a$ is a filling
of the $n$-boundary $b$ if $\partial a = b$.  The filling length of $b$ is the least number of cells needed
to fill it.  Denote the filling length of $b$ by $\ell_f(b)$, and the number of cells in $b$ by $|b|$.  
The $n$-th Dehn function of $X$ is the function $d^n : \N \rightarrow \R_+$ defined by
\[ d^n( k ) = \max_{b} \ell_f(b) \]
where this maximum is taken over $n$-boundaries $b$ with $|b| \leq k$.  In this way given an $n$-boundary $b$
it is possible to find a filling $a$ with $|a| \leq d^n( | b | )$.  This notion does not take into account
the position of $b$ in $X$, only how many cells in $b$.  

Let $x_0$ be a fixed vertex of $X$.  This induces a length function on the vertices by $\ell_X(v) = d_X(x_0,v)$,
where $d_X$ is the graph metric on $X^{(1)}$.  Let $\sigma$ be an $n$-cell with vertices $v_0$, $v_1$, \ldots, $v_n$.
define the length of $\sigma$ to be $\ell_X( \sigma ) = \ell_X( v_0 ) + \ldots + \ell_X( v_n )$, the sum
of the lengths of the vertices.  The weighted number of cells in $b$ is given by
\[ |b|_w = \sum_{\sigma \in b} \ell_X( \sigma ) \]
For a boundary $b$ the weighted filling length of $b$, $\ell^w_f(b)$, is $\min \{ | a |_w \, | \, \partial a = b \}$.
The $n$-th weighted Dehn function of $X$, $d_w^n : \N \rightarrow \R_+$, is given by
\[ d_w^n( k ) = \max \left\{ \ell^w_f(b) \, | \, |b|_w \leq k \right\} \]
In this way if $b$ is an $n$-boundary, there is a filling $a$ with $|a|_w \leq d_w^n( |b|_w )$.
A function $f: \N \rightarrow \R_+$ dominates $g : \N \rightarrow \R_+$ if there are constants $A$, $B$, $C$, $D$, and $E$
such that for all $n$, $f(n) \leq A g( B n + C ) + Dn + E$.  $f$ and $g$ are equivalent if each dominates the other.
In the case of the usual Dehn function, this is the natural notion of equivalence.

\begin{lem}
Suppose that $X$ and $Y$ are two polyhedral complexes with connected $1$-skeletons, upon which a finitely generated
group $G$ acts properly by isometries, each with only finitely many orbits of cells in every dimension. 
Then $X$ and $Y$ have equivalent weighted Dehn functions in all dimensions.  
\end{lem}
\proof
Let the fixed base-point vertices of $X$ and $Y$ be $x_0$ and $y_0$ respectively.
$X^{(1)}$ and $Y^{(1)}$ are quasiisometric, as they are each quasiisometric to $G$.
Let $\Phi : X^{(1)} \rightarrow Y^{(1)}$ and $\Psi : Y^{(1)} \rightarrow X^{(1)}$ be
quasiinverse quasiisomorphisms.  At the expense of enlarging the constants involved, it is
assumed that $\Phi$ and $\Psi$ map vertices to vertices and basepoint to basepoint.  Let $d_{w,X}^1$ 
and $d_{w,Y}^1$ be the weighted Dehn functions for filling $1$-boundaries in $X$ and $Y$ respectively.  As 
there are only finitely many orbits of $2$-cells in $X$ and $Y$, assume that no $2$-cell 
has more than $J$ vertices.

Let $y_1, y_2, \ldots, y_{n}$ be the vertices of some $1$-boundary, $\beta$, in $Y^{(1)}$.  Let $v_i = \Psi( y_i )$.
There is a constant $C$ such that for any $y, y' \in Y^{(0)}$, $d_X( \Psi(y), \Psi(y') ) \leq C d_Y( y, y')$.
As $d_Y(y_i, y_{i+1}) = 1$ and $d_Y(y_n, y_1) = 1$, $d_X( v_i, v_{i+1} ) \leq C$  and $d_X( v_n, v_1 ) \leq C$.
Let $\alpha_{i,i+1}$ be a geodesic path in $X^{(1)}$ connecting $v_i$ to $v_{i+1}$, and $\alpha_{n,1}$
connecting $v_n$ to $v_1$.  The concatenation of the $\alpha_{i,j}$ paths yield a $1$-boundary, $\alpha$,
in $X^{(1)}$.  Each vertex in one of the $\alpha_{i,i+1}$ paths is within $C$ from $v_i$, so has
length no more than $\ell_X(v_i) + C$.  The weighted length of each edge in $\alpha_{i,i+1}$
is no more than $2\ell_X(v_i) + C$.  As there are at most $C$ edges in each $\alpha_{i,i+1}$,
it has weighted length bounded by $2C \ell_X(v_i) + C$.  The weighted length of $\alpha$ is 
bounded by $\sum_{i} \left(2C \ell_X(v_i) + C \right)$.  If $\lambda$ and $D$ are the quasiisometry
constants of $\Psi$, then $d_X( x_0, v_i ) \leq \lambda d_Y( y_0, y_i ) + D$, so 
the weighted length of $\alpha$ is bounded by $\sum_{i} \left(2C \lambda \ell_Y(y_i) + C(2D+1) \right)$.
As the weighted length of $\beta$ is $2 \sum_i \ell_Y(y_i)$ there is a positive constant $M$ 
such that $| \alpha |_w \leq M | \beta |_w$, with $M$ independent of $\beta$.

There is a filling $\gamma$ with $| \gamma |_w \leq d_{w,X}^1( |\alpha|_w ) \leq d_{w,X}^1( M |\beta|_w )$.
Let $\sigma$ be a $2$-cell in $\gamma$, with vertices $(x_1, x_2, \ldots, x_j)$.  Let $u_i = \Phi( x_i )$.
As $d_X(x_i, x_{i+1}) = 1$ and $d_X( x_j, x_1 ) = 1$, $d_Y( u_i, u_{i+1} ) \leq C'$ and $d_Y( u_n, u_1 ) \leq C'$
for some constant $C'$ as above.  Let $\mu_{i,i+1}$ and $\mu_{j,1}$ be geodesic paths in $Y^{(1)}$ connecting
$u_i$ to $u_{i+1}$, and connecting $u_j$ to $u_1$, respectively.  Denote the resulting cycle obtained by concatenating
these paths, by $\mu$.  $\mu$ need not be a $2$-cell in $Y$, however the number of vertices in the boundary of 
$\sigma$ is bounded by $J$.  Thus the length around the boundary of $\mu$ is bounded by at most $C' J$.
There is a constant $L$ such that each $\mu$ can be filled by at most $L$ $2$-cells in $Y$. Denote a filling with
minimal weighted length as $\Phi( \sigma )$.  $\Phi( \sigma )$ is not a single $2$-cell, but is a connected
subcomplex of $Y^{(2)}$.  The length of each vertex in $\Phi( \sigma )$ is bounded by $\ell_Y( u_1 ) + JL$.
As above we find a positive constant $M'$ such that $|\Phi( \sigma )|_w \leq M' | \sigma |_w$.
Let $\Phi( \gamma )$ be the subcomplex spanned by all of the $\Phi( \sigma )$, for all $2$-cells $\sigma \in \gamma$.
$| \Phi( \gamma ) |_w \leq M' | \gamma |_w$, with $\Phi( \gamma )$ nearly filling $\beta$.

Consider $\Phi( v_i) = \Phi( \Psi( y_i ) )$.  There is $K$ such that $d_Y( y_i, \Phi(v_i) ) \leq K$ for all $i$.
Let $\nu_i$ be a geodesic path in $Y$ connecting $y_i$ to $\Phi( v_i )$.  Recall that $\alpha_{i,i+1}$ connects
$v_i$ to $v_{i+1}$ in $X$.  In the construction of $\Phi( \gamma )$, each edge of $\alpha_{i,i+1}$ was lifted
back into $Y$, as a path of edges, say $\Phi( \alpha_{i,i+1} )$.  The concatenation of the edge between 
$y_i$ and $y_{i+1}$, $[y_i, y_{i+1}],$ with $\nu_i$, $\Phi( \alpha_{i,i+1} )$, and $\nu_{i+1}$ 
gives a $1$-boundary in $Y$ with uniformly bounded length around the cycle.  There is then a filling of this 
cycle, $\eta_{i,i+1}$, by a uniformly bounded number  of $2$-cells in $Y$, and 
$| \eta_{i,i+1} |_w \leq M'' | [y_i, y_{i+1}] |_w$, for some universal constant $M''$.  
Let $\beta_f$ be the subcomplex spanned by $\Phi( \gamma )$ and all of the $\eta_{i,i+1}$.
From the construction, $\beta_f$ is a filling of $\beta$, and 
$| \beta_f |_w \leq | \Phi( \gamma ) |_w + \sum | \eta_{i,i+1} |_w$.
Also from above, $| \beta_f |_w \leq M' d_{w,X}^1( M | \beta |_w ) + M'' | \beta |_w$.
It follows that a weighted Dehn function for $X$ dominates one for $Y$.  By reversing the quasiisometries
the two $1$-dimensional weighted Dehn functions are seen to be equivalent.  The higher dimensional cases
are proven similarly.
\endproof

By examining the proof, replacing the weighted lengths by the usual cardinality counting the following is also apparent.
\begin{lem}
Suppose that $X$ and $Y$ are two polyhedral complexes acted upon properly by a discrete group $G$.
Moreover assume that in each dimension, $X$ and $Y$ have only finitely many orbits of cells under this action,
and have connected $1$-skeletons.  $X$ and $Y$ have equivalent Dehn functions in all dimensions.
\end{lem}
For the case of filling $1$-boundaries, this is well-known.  \cite{BRS}

A natural question is how a weighted Dehn function relates to the usual Dehn function.
\begin{lem}\label{lemBddWeightedDehn}
Suppose that $G$ and $X$ are as in the previous lemma.  $d_w^n(x)$ is bounded above by 
$d^n(x) \left( x + d^n(x) \right)$, up to equivalence.
\end{lem}
\proof
Let $u$ be an $n$-boundary in $X$, with weighted length $|u|_w$.  There is a filling $\omega$ of $u$ by $(n+1)$-cells
with $N = | \omega | \leq d^n( |u| ) \leq d^n( |u|_w )$ as above.  We estimate $| \omega |_w$ in terms of $| u |_w$.
Let $\sigma_1$, $\sigma_2$, $\ldots$, $\sigma_N$ be the $(n+1)$-cells of $\omega$.  By the finiteness property of
$X$, there are constants $J$ and $J'$ such that each $\sigma_i$ has no more than $J$ edges and $J'$ vertices.  
Let $v$ be a vertex of $u$.  For every vertex $v'$ of any $\sigma_i$, $\ell_X( v' ) \leq \ell_X(v) + NJ$.
\begin{eqnarray*}
| \sigma_i |_w & = & \sum_{v' \in \sigma_i} \ell_X( v' )\\
 & \leq & \sum_{v' \in \sigma_i} \left(\ell_X( v ) + NJ \right)\\
 & \leq & J' \ell_X( v ) + NJJ'\\
 & \leq & J' | u |_w + JJ' d^n( |u| )\\
 & \leq & J' | u |_w + JJ' d^n( |u|_w )
\end{eqnarray*}
As there are no more than $d^n( |u|_w )$ such $\sigma_i$, we have that
$d_w^n( |u|_w ) \leq d^n( |u|_w ) \left( J' |u|_w + JJ' d^n( |u|_w ) \right)$
\endproof

\begin{lem}
Let $G$ and $X$ be as above.  $d^n(x)$ is bounded by $x d_w^n( x(x+1) )$, up to equivalence.
\end{lem}
\proof
Let $x_0$ be the basepoint of $X$, and let $u$ be a connected $n$-boundary in $X$.
By cocompactness there is a positive constant $L$ such that for any vertex $v \in X$ there is $g \in G$ such that  
$d_X( x_0, g v ) \leq L$.  Let $v$ be a vertex of $u$, let $g$ be such a group element, and let $u' = g \cdot u$.
There are constants $J$ and $J'$ such that any $n$-cell of $X$ has at most $J$ edges and $J'$ vertices.
The length of each vertex of $u'$ is bounded by $L + J |u|$. The weighted length of each cell of $u'$ is
then bounded by $J' \left( L + J |u| \right)$, yielding $| u' |_w \leq J' |u| \left( L + J |u| \right)$.
Let $\omega'$ be a filling of $u'$ with $| \omega' |_w \leq d_w^n( |u'|_w )$.  $\omega = g^{-1} \cdot \omega'$
is a filling of $u$ with 
\begin{eqnarray*}
| \omega | & = & | \omega' | \\
 & \leq & | \omega' |_w \\
 & \leq & d_w^n( |u'|_w ) \\
 & \leq & d_w^n\left( J' |u| \left( L + J |u| \right) \right)
\end{eqnarray*}

Assume that $u$ is not connected.  Let $u_1$, $u_2$, $\ldots$, $u_k$ be the connected components of $u$.
Each $u_i$ is itself a boundary so there are fillings $\omega_i$ of $u_i$ with
$| \omega_i | \leq d_w^n\left( J' |u_i| \left(L + J |u_i| \right) \right)$.  As $|u_i| \leq |u|$ and $k \leq |u|$
the collection of all $\omega_i$ constitute a filling of $u$ with no more than 
$|u| d_w^n\left( J' |u| \left( L + J |u| \right) \right)$ cells.
\endproof

\begin{cor}
Let $G$ and $X$ be as above.  The Dehn function $d^n$ is polynomially bounded if and only if the weighted 
Dehn function $d_w^n$ is polynomially bounded.
\end{cor}
This shows that there is a wide class of groups having polynomially bounded weighted Dehn functions in all dimensions.  
By work of Gersten it includes all finitely generated groups endowed with combings of polynomial length. 
We will also show later that groups of polynomial growth are in this class.

Suppose $G$ has a classifying space $X'$ with the homotopy type of a polyhedral complex with
finitely many cells in each dimension, and let $X$ be the universal cover of $X'$.  
This is the case for combable groups by the work of Gersten in \cite{ger}.
Ogle refers to such groups as $HF^\infty$ groups in \cite{Ogle}.  The following characterizes the isocohomological 
$HF^\infty$ groups, in terms of higher Dehn functions.

\begin{thm}\label{thmClassification}
For an $HF^\infty$ group $G$, with $X$ as above, the following are equivalent.
\begin{enumerate}
\item[(1)] $HP^*(G; V) \rightarrow H^*(G; V)$ is an isomorphism for all coefficients $V$.
\item[(2)] $HP^*(G; V) \rightarrow H^*(G; V)$ is surjective for all coefficients $V$.
\item[(3)] All higher Dehn functions of $X$ are polynomially bounded.
\end{enumerate}
\end{thm}

\proof
(1) $\implies$ (2): Obvious.  

(2) $\implies$ (3):  This is similar in spirit to Mineyev's proof that if the comparison homomorphism from
bounded cohomology of $G$ to the group cohomology of $G$ is surjective for all coefficients, then $G$ is 
hyperbolic \cite{min}.  The difficulty here lies in analyzing the Fr\'echet nature of our coefficients, rather than
the Banach structure in the bounded case.

Let $V$ be the collection of polynomially bounded $(n-1)$-boundaries with complex coefficients, endowed with 
the Fr\'echet space structure endowed by the family of filling norms defined as follows:
Let $\xi \in V$.  We say $\| \xi \|_{f,k} \leq M_k$ if there is a $\phi \in C_n(X)$ of polynomial growth
with $\partial \phi = \xi$, and $\| \phi \|_k \leq M_k$.  

Let $Y$ be the geometric realization of the bar complex of $G$.  That is, there is an $n$-cell in $Y$ for every
$(n+1)$-tuple $[g_0, g_1, \ldots, g_n]$ of elements of $G$, endowed with the diagonal action of $G$.  
Let $C_n(X)$ and $C_n(Y)$ denote the complex valued algebraic
$n$-chains in $X$ and $Y$ respectively.  As both complexes yield projective $\CG$-module resolutions of $\C$, 
there are homotopy equivalences $\phi_* : C_*(Y) \rightarrow C_*(X)$ and $\psi_* : C_*(X) \rightarrow C_*(Y)$
which are $\CG$-module morphisms.

Let $C^*(X,V)=\Hom_{\CG}( C_*(X), V )$ and $C^*(Y,V)=\Hom_{\CG}( C_*(Y), V )$ be the dual cochain complexes
with dual maps $\phi^* : C^*(X,V) \rightarrow C^*(Y,V)$ and $\psi^* : C^*(Y,V) \rightarrow C^*(X,V)$.
The cochain map $\psi^* \circ \phi^*$ is homotopic to the identity, so $\psi^* \circ \phi^*$ induces the identity
map on cohomology $H^*( G, V )$ in all positive dimensions.  We will make use of the pairings 
$< \cdot | \cdot > : C^*(X,V) \oplus C_*(X) \rightarrow V$ and   
$< \cdot | \cdot > : C^*(Y,V) \oplus C_*(Y) \rightarrow V$.

Consider the map $u : C_n(X) \rightarrow V$ given by the composition 
$C_n(X) \stackrel{\partial}{\rightarrow} B_{n-1}(X) \hookrightarrow V$,
where $B_{n-1}(X)$ is the image of $\partial : C_n(X) \rightarrow C_{n-1}(X)$; They are the finitely 
supported $(n-1)$-boundaries, which can be included into $V$.  $u$ is a linear map, commuting with the 
$\CG$-action, so $u \in C^n(X,V)$.  In fact $u$ is an $n$-cocycle.
As $\psi^* \circ \phi^*$ is the identity map in cohomology, there is a $(n-1)$-cochain $v$ with 
$u = (\psi^n \circ \phi^n)(u) + \delta v$.
$\phi^n(u)$ is a $n$-cocycle in $C^n(Y,V)$, so by assumption there is a polynomially bounded 
$n$-cocycle $u'$ and some $(n-1)$-cochain $v'$
such that $\phi^n(u) = u' + \delta v'$.  As $V$ is a Fr\'echet space with a family of norms 
$\| \cdot \|_{f,k}$, we must make precise
exactly what is meant that $u'$ is polynomially bounded.  It means the following:  
For every $k$ there exists a polynomial 
$P_k$ such that for each $\xi \in C_n(Y)$, 
$\| u'( \xi ) \|_{f,k} \leq  P_k( \| \xi \|_{k} )$, where $\| \xi \|_{k}$ is the usual polynomially 
weighted $\ell^1$ norm.  

For an $(n-1)$-chain 
\[b = \sum_{g_0,g_1, ..., g_{n-1} \in G} \beta_{g_0,g_1, ..., g_{n-1}} [ g_0, g_1, ..., g_{n-1}]\]
in $Y$, define the cone over $b$, $[e,b]$, to be
\[ [e,b] = \sum_{g_0,g_1, ..., g_{n-1} \in G} \beta_{g_0,g_1, ..., g_{n-1}} [ e, g_0,g_1, ..., g_{n-1} ] \]
It is clear that if $b$ is a cycle, $\partial [e,b] = b$ so $[e,b]$ serves as a filling of $b$, with the 
property that $\| [e,b] \|_k = \| b \|_k$ for all $k$.
If $\alpha$ is a cocycle in $C^n(X,V)$ and $c \in C_n(X)$ then $<\alpha \, | \, c - [e,\partial c]> = 0$, since 
$c - [e,\partial c]$ is a boundary.  Therefore $<\alpha \, | \, c > = < \alpha \, | \, [e,\partial c] >$.

Let $b$ be any $(n-1)$-boundary in $B_{n-1}(X)$, and let $a$ be any $n$-chain with $\partial a = b$.
\begin{eqnarray*}
b & = & \partial a \\
  & = & < u \, | \, a >\\
  & = & < (\psi^n \circ \phi^n)(u) + \partial v \, | \, a >\\
  & = & < (\psi^n \circ \phi^n)(u) \, | \, a > + < v \, | \, b >\\
  & = & < \phi^n(u) \, | \, \psi_n(a) > + < v \, | \, b >\\
  & = & < \phi^n(u) \, | \, [ e, \partial \psi_n(a)] > + < v \, | \, b >\\
  & = & < \phi^n(u) \, | \, [ e, \psi_{n-1}(b)] > + < v \, | \, b >\\
  & = & < u' + \delta v' \, | \, [ e, \psi_{n-1}(b)] > + < v \, | \, b >\\
  & = & < u' \, | \, [ e, \psi_{n-1}(b)] > + < v' \, | \, \partial [ e, \psi_{n-1}(b)] > + < v \, | \, b >\\
  & = & < u' \, | \, [ e, \psi_{n-1}(b)] > + < v' \, | \, \psi_{n-1}(b) > + < v \, | \, b >\\
  & = & < u' \, | \, [ e, \psi_{n-1}(b)] > + < \psi^{n-1}(v') \, | \, b > + < v \, | \, b >\\
  & = & < u' \, | \, [ e, \psi_{n-1}(b)] > + < \psi^{n-1}(v') + v \, | \, b >
\end{eqnarray*}

Thus for each $k$,
\begin{eqnarray*}
\| b \|_{f,k} & \leq & \| < u' \, | \, [ e, \psi_{n-1}(b)] > \|_{f,k} + \| < \psi^{n-1}(v') + v \, | \, b >\|_{f,k}\\
 & \leq & P_k( \| \psi_{n-1}(b) \|_{k} ) + \| < \psi^{n-1}(v') + v \, | \, b >\|_{f,k}
\end{eqnarray*}

$\psi_{n-1} : C_{n-1}(X) \rightarrow C_{n-1}(Y)$ is a linear map, commuting with the
$G$-action.  As $X$ has finitely many orbits of $(n-1)$-cells, there are 
$\sigma_0$, $\ldots$ , $\sigma_i \in X^{(n-1)}$, representatives
of the orbits, such that for each $\sigma \in G \cdot \sigma_j$, $\ell_X( \sigma_j ) \leq \ell_X( \sigma )$.
Let $g \in G$ be such that $\sigma = g \cdot \sigma_j$, with $\sigma_j$ having vertices $( v_0, \ldots, v_{n-1} )$.
\begin{eqnarray*}
\ell_X( \sigma ) & = & \ell_X(g \sigma_j) \\
 & = & \ell_X( g v_0 ) + \ldots + \ell_X( g v_{n-1} )\\
 & = & d_X( x_0, g v_0 ) + \ldots + d_X( x_0, g v_{n-1} )\\
 & \leq & d_X( x_0, g x_0 ) + d_X(g x_0, g v_0 ) + \ldots \\
 & & + d_X( x_0, g x_0 ) + d_X( g x_0, g v_{n-1} )\\
 & = & n d_X( x_0, g x_0 ) + d_X( x_0, v_0 ) + \ldots + d_X( x_0, v_{n-1}\\
 & = & n d_X( x_0, g x_0 ) + \ell_X( \sigma_j )
\end{eqnarray*}
Since the quotient of $X^{(1)}$ by the $G$-action is finite, the length function on $G$ given by
$\hat{\ell}(g) = d_X( x_0, g x_0 )$ is bilipschitz equivalent to the word-length function on $G$.
Thus there is a constant $C' > 0$ such that $\frac{1}{C'} \ell(g) \leq d_X( x_0, g x_0 ) \leq C' \ell( g )$.
In this way, $\ell_X( \sigma ) \leq n C' \ell(g) + \ell_X( \sigma_j )$.  As there
are only finitely many $\sigma_j$, there are constants $C$ and $D$ such
that $\ell_X( \sigma ) \leq C \ell(g) + D$.

Similarly, 
\begin{eqnarray*}
\ell_X( \sigma ) & = & \ell_X(g \sigma_j) \\
 & = & \ell_X( g v_0 ) + \ldots + \ell_X( g v_{n-1} )\\
 & = & d_X( x_0, g v_0 ) + \ldots + d_X( x_0, g v_{n-1} )\\
 & \geq & d_X( x_0, g x_0 ) - d_X(g x_0, g v_0 ) + \ldots \\
 &  & + d_X( x_0, g x_0 ) - d_X( g x_0, g v_{n-1} )\\
 & \geq & n d_X( x_0, g x_0 ) - \left( d_X( x_0, v_0 )  + \ldots + d_X( x_0, v_{n-1} ) \right)\\
 & \geq & n d_X( x_0, g x_0 ) - \ell_X( \sigma_j )
\end{eqnarray*}
Combining this with the bilipschitz constant from above yields
\begin{eqnarray*}
\ell(g) & \leq & C' d_X( x_0, g x_0 ) \\
 & \leq & \frac{C'}{n} \ell_X( \sigma) + \frac{1}{n} \ell_X( \sigma_j )
\end{eqnarray*}
As there are only finitely many $\sigma_j$, there are constants $A$ and $B$ such that
$\ell(g) \leq A \ell_X( \sigma ) + B$.
Thus the lengths $\ell(g)$ and $\ell_X( \sigma )$ are linearly equivalent.

Again consider $\| \psi_{n-1}(b) \|_k$, with 
\[ \psi_{n-1}( \sigma_j ) = \sum_{ g_0, \ldots, g_{n-1} } \alpha^j_{ g_0, \ldots, g_{n-1} } [ g_0, \ldots, g_{n-1} ]\]
Then $\psi_{n-1}( g \cdot\sigma_j ) = \sum_{ g_0, \ldots, g_{n-1} } \alpha^j_{ g_0, \ldots, g_{n-1} } [ g g_0, \ldots, g g_{n-1} ]$
by the equivariance of $\psi_{n-1}$.  In particular, 
\begin{eqnarray*}
\| \psi_{n-1}( g \cdot \sigma_j ) \|_k & = & \sum_{ g_0, \ldots, g_{n-1} } |\alpha^j_{ g_0, \ldots, g_{n-1} }| 
\left( 1 + \ell( g g_0 ) + \ldots + \ell( g g_{n-1} ) \right)^k\\
 & \leq & \sum_{ g_0, \ldots, g_{n-1} } |\alpha^j_{ g_0, \ldots, g_{n-1} }| 
\left( 1 + \ell( g_0 ) + \ldots + \ell( g_{n-1} ) + n \ell(g) \right)^k\\
& \leq & \sum_{ g_0, \ldots, g_{n-1} } |\alpha^j_{ g_0, \ldots, g_{n-1} }| 
\left( 1 + \ell( g_0 ) + \ldots + \ell( g_{n-1} ) \right)^k \cdot n^k \left( 1 + \ell(g) \right)^k\\
& = & n^k \left( 1 + \ell(g) \right)^k \| \phi_{n-1}(\sigma_j) \|_k
\end{eqnarray*}
For each $k$, let $\| \phi_{n-1} \|_{\infty,k} = \max \left\{ \| \phi_{n-1}( \sigma_j ) \|_k \, | \, j = 0, 1, \ldots, i \right\}$, a finite constant.
Then if $b = \sum_{\sigma \in X^{(n-1)}} \beta_\sigma \cdot \sigma$, let $b_j = \sum_{g \in G} \beta_{g, \sigma_j} \cdot g \sigma_j$ be the projection
onto the orbit spanned by $\sigma_j$.  As there may be several group elements translating $\sigma_j$ to $\sigma$, we pick one of minimal length
and let $\beta_{g, \sigma_j} = \beta_\sigma$ for this particular $g$.  For all other group elements translating $\sigma_j$ to $\sigma$, we require
$\beta_{g, \sigma_j}=0$.  In this way there is a single group element representing every cell in this orbit, and $\| b \|_k = \sum_j \| b_j \|_k$.
\begin{eqnarray*}
\| \psi_{n-1}( b_j ) \|_k & = & \| \sum_{g \in G} \beta_{g, \sigma_j} \psi_{n-1}( g \cdot \sigma_j ) \|_k \\
 & \leq & \sum_{g \in G} | \beta_{g, \sigma_j} | \| \psi_{n-1}( g \cdot \sigma_j ) \|_k\\
 & \leq & n^k \| \psi_{n-1} \|_{\infty,k} \sum_{g \in G} | \beta_{g, \sigma_j} |\left( 1 + \ell(g) \right)^k \\
 & \leq A \| b_j \|_k 
\end{eqnarray*}
As $\psi_{n-1}(b) = \sum_j \psi_{n-1}( b_j )$, we have that there is a constant $A_k$ such that 
$\| \psi_{n-1}( b ) \|_k \leq A_k \| \psi_{n-1} \|_{\infty,k} \| b \|_k$, whence the first term in the above is appropriately bounded.

Now consider a linear $G$-equivariant map $\eta :  C_{n-1}(X) \rightarrow V$.  We wish to bound $\| \eta( b ) \|_{f,k}$
in terms of $\| b \|_k$.  As above, let $b = \sum_j b_j$.
Let $\| \eta \|_{f,\infty,k} = \max \left\{ \| \eta( \sigma_j ) \|_{f,k} \right\}$, and $a_j \in C_n(X)$ be such
that $\| \eta( \sigma_j ) \|_{f,k} = \| a_j \|_k$ with $\partial a_j = \eta( \sigma_j )$.
Then $\partial( g \cdot a_j ) = \eta( g \cdot \sigma_j )$.  Let 
$a_j = \sum_{\gamma \in X^{(n-1)}} \alpha^j_{\gamma} \gamma$.
\begin{eqnarray*}
\| \eta( g \cdot \sigma_j ) \|_{f,k} & \leq & \| g \cdot a_j \|_k \\
 & = & \| \sum_{\gamma \in X^{(n-1)}} \alpha^j_{\gamma} g \cdot\gamma \|_k \\
 & = & \sum_{\gamma \in X^{(n-1)}} |\alpha^j_{\gamma}| \left( 1 + \ell_X( g \cdot\gamma ) \right)^k\\
 & \leq & \sum_{\gamma \in X^{(n-1)}} |\alpha^j_{\gamma}| 
	\left( 1 + \ell_X( g v^\gamma_0 ) + \ldots +  \ell_X( g v^\gamma_{n-1} )\right)^k\\
 & \leq & M \left( 1 + \ell(g) \right)^k \sum_{\gamma \in X^{(n-1)}} |\alpha^j_{\gamma}| 
	\left( 1 + \ell_X( v^\gamma_0 ) + \ldots +  \ell_X( v^\gamma_{n-1} )\right)^k\\
 & = & M \left( 1 + \ell(g) \right)^k \| a_j \|_k \\
 & = & M \left( 1 + \ell(g) \right)^k \| \eta( \sigma_j ) \|_{f,k}\\
 & \leq & M \left( 1 + \ell(g) \right)^k \| \eta \|_{f,\infty,k}
\end{eqnarray*}
for some positive constant $M$.  Let $b_j$ be as above.
\begin{eqnarray*}
\| \eta( b_j ) \|_{f,k} & = & \| \sum_{g \in G} \beta_{g, \sigma_j} \eta( g \cdot \sigma_j ) \|_{f,k} \\
 & \leq & \sum_{g \in G} | \beta_{g, \sigma_j} \| \eta( g \cdot \sigma_j ) \|_{f,k}\\
 & \leq & M \| \eta \|_{f,\infty,k} \sum_{g \in G} | \beta_{g, \sigma_j} | \left( 1 + \ell(g) \right)^k\\
 & \leq & M' \| \eta \|_{f,\infty,k} \| b_j \|_k
\end{eqnarray*}
As before it follows that there is a constant $B_k$ such that $\| \eta( b ) \|_{f,k} \leq B_k \| b \|_k$.
Letting $\eta = \psi^{n-1}(v') + v$ the second term of the above is appropriately bounded.
This proves the result in filling from dimension ${n-1}$ to 
dimension $n$.

(3) $\implies$ (1):
This implication is similar to Theorem 2.2.4 of Ogle in \cite{Ogle}.  
Ogle's definition of higher Dehn functions is somewhat different than that given here, so our 
proofs are different.  We work in the framework of bornologies as developed in \cite{Hn, meyer3}.
Let $\SnX[n]$ be the set of all function $\phi : X^{(n)} \rightarrow \C$
for which each of the norms 
\[ \| \phi \|^X_k = \sum_{ \sigma \in X^{(n)}} | \phi( \sigma ) | \left( 1 + \ell_X( \sigma ) \right)^k \]
is finite.   This is the $\ell^1$ rapid-decay completion of the space of complex-valued $n$-chains, $\CnX[n]$.
We endow $\SnX[n]$ with the Fr\'echet bornology induced by this family of norms.
The $G$ action on $X^{(n)}$ extends to a boronological $\SG$-module structure on $\SnX[n]$.  With the
usual boundary map we find a complex of bornological $\SG$-modules:
\[ \ldots \rightarrow \SnX[3] \rightarrow \SnX[2] \rightarrow \SnX[1] \rightarrow \SnX[0] \rightarrow \C \rightarrow 0 \]The minimal fillings allow us to map $n$-boundaries into $\SnX[n+1]$, due to the polynomially bounded Dehn 
functions.  Diagram chasing extends this map to a bounded $\SnX[n] \rightarrow \SnX[n+1]$, yielding a $\C$-splitting
of the complex.  
In particular, using the minimal fillings we construct a bounded linear map 
$f_n : \ker \partial \rightarrow \CnX[n+1]$
for which if $\xi \in \CnX[n]$ is in the image of $\partial : \CnX[n+1] \rightarrow \CnX[n]$ then
$\partial f_n( \xi ) = \xi$.  Given an arbitrary $\xi \in \CnX[n]$, $\xi$ may not be in $\ker \partial$, but
$\partial \xi$ is.  Consider the map
$s_n : \CnX[n] \rightarrow \CnX[n+1]$ by the formula $s_n( \xi ) = f_n\left( \xi - f_{n-1}( \partial \xi ) \right)$,
using $\partial \left( \xi - f_{n-1}( \partial \xi ) \right) = 0$.
In this way
\begin{eqnarray*}
(s_{n-1} \partial + \partial s_n)(\xi) & = & s_n( \partial \xi ) + \partial s_n( \xi )\\
 & = & f_{n-1} \left( \partial \xi - f_{n-2}(\partial \partial \xi ) \right) + 
	\partial \left( f_n\left( \xi - f_{n-1}( \partial \xi ) \right) \right)\\
 & = & f_{n-1}( \partial \xi ) + \left( \xi - f_{n-1}) \partial \xi \right)\\
 & = & \xi
\end{eqnarray*}
As the $s_n$ maps are linear, being the composition of linear maps, they form a $\C$-splitting
of the $\CnX[*]$ complex.  The polynomially bounded weighted Dehn function ensures that this splitting is bounded
in the bornology on $\CnX[n]$ induced as a subspace of $\SnX[n]$, so it extends to a bounded $\C$-linear
splitting of the $\SnX[*]$ complex.
Compare this complex to the bornologically projective $\CG$-resolution:
\[ \ldots \rightarrow \CnX[3] \rightarrow \CnX[2] \rightarrow \CnX[1] \rightarrow \CnX[0] \rightarrow \C \rightarrow 0 \]
where each of the $\CnX[n]$ are endowed with the fine bornology.

As there are only finitely many $G$-orbits of $n$-cells, $\CnX[n]$ is bornologically isomorphic to $\CG \btensor W$ as 
bornological $\CG$-modules, for $W$ a finite dimensional $\C$-vector space, spanned by the $G$-orbits of $n$-cells.
Similarly $\SnX[n] \isom \SG \btensor W$  as bornological $\SG$-modules,for the same space $W$, implying that
$\SnX[*]$ is a bornologically projective resolution of $\C$ over $\SG$.  

Let $V$ be any bornological $\SG$-module.  $V$ is then a bornological $\CG$-module by restriction.
Using the properties of the bounded homomorphism functor we find
\begin{eqnarray*}
\bHom_{\SG}( \SnX[n], V ) & \isom & \bHom_{\SG}( \SG \btensor W, V )\\
 & \isom & \bHom( W, V )\\
 & \isom & \bHom_{\CG}( \CG \btensor W, V )\\
 & \isom & \bHom_{\CG}( \CnX[n], V )
\end{eqnarray*}
As $\bHom_{\SG}( \SnX[n], V)$ and $\bHom_{\CG}( \CnX[n], V )$ are isomorphic, the cohomology
of the cocomplexes obtained by applying $\bHom_{\SG}( \cdot, V )$ and $\bHom_{\CG}( \cdot, V )$
are equal.  By the work of Meyer, these are the polynomial cohomology of $G$ and the group cohomology of $G$, respectively.
\endproof

By \cite{ger2} finitely generated nilpotent groups have polynomially bounded Dehn function for filling $1$-boundaries.
In fact, nilpotent groups are isocohomological \cite{Ji, Ogle, meyer2} and are of type $HF^\infty$ so we have 
established the following generalization.
\begin{cor}
Finitely generated nilpotent groups have all higher Dehn functions polynomially bounded.
\end{cor}

We have also established the following.
\begin{cor}
Let $\Z$ act on $\Z^2$ by the matrix 
$\begin{pmatrix}
  2 & 1\\
  1 & 1
\end{pmatrix}$.  $\Z^2 \rtimes \Z$ is not isocohomological.
\end{cor}
\proof
It is clear that the semidirect product is of type $HF^\infty$.  By 
\cite{Epstein} this group has Dehn function which is at least exponential.
\endproof
This is the first example of a group without the isocohomological property.

The following is a natural generalization of the fact that a group quasiisometric to a
finitely presented group is itself finitely presented \cite{GH}.
\begin{lem}\label{lemHFinfty}
The class of $HF^\infty$ groups is closed under quasiisometry.
\end{lem}
\proof
Let $G$ and $H$ be two quasiisometric groups, with $G$ of type $HF^\infty$.  Let $X$ be a polyhedral complex on 
which $G$ acts freely by isometries, with finitely many orbits of cells in each dimension, as guaranteed by the
$HF^\infty$ property.  Let $\Phi : G \rightarrow H$ and $\Psi : H \rightarrow G$ be the coarse inverse 
quasiisometries.  $G$ is also quasiisometric to $X^{(1)}$.  Denote by $\alpha : G \rightarrow X^{(1)}$
and $\beta : X^{(1)} \rightarrow G$ the coarse inverse quasiisometries.  The compositions $\alpha \circ \Psi$
and $\Phi \circ \beta$ are coarse inverse quasiisometries between $H$ and $X^{(1)}$.  Let $\lambda$ and $C$
be quasiisometry constants for each of these maps and 
$\left( \Phi \circ \beta \right) \circ \left( \alpha \circ \Psi \right)$.

We build an acyclic polyhedral complex on which $H$ acts freely by isometries, with finitely many orbits of 
cells in each dimension.  Let $\Gamma^{(1)}$ be the Cayley graph of $H$ with respect to some fixed generating set.
We identify the elements of $H$ with the corresponding vertices in $\Gamma^{(0)}$.

There is a positive constant $L_2$ such that every elementary $2$-cell in $X^{(2)}$ has at most $L_2$ edges around
its boundary.  Let $u$ be such a $2$-cell, and denote the vertices, around its boundary as 
$x_1$, $x_2$, $\ldots$, $x_n$, $x_{n+1}=x_1$, with $d_{X^{(1)}}(x_i, x_{i+1}) = 1$.  Let $x'_i = \Phi \circ \beta(x_i)$.
$d_{H}( x'_i, x'_{i+1} ) \leq \lambda + C$.  Connecting $x'_i$ to $x'_{i+1}$ with a geodesic in $\Gamma^{(1)}$
we obtain a cycle of edges of length at most $L_2 \left(\lambda + C\right)$.  As $\alpha \circ \Psi$ and
$\Phi \circ \beta$ are coarse inverses, there is a $K$ such that for all $h \in H$, 
$d_{H}(h, \left( \Phi \circ \beta \right) \circ \left( \alpha \circ \Psi \right)(h) ) \leq K$.
If $d_H( h,h' ) = 1$ then 
$d_H( \left( \Phi \circ \beta \right) \circ \left( \alpha \circ \Psi \right)(h), \left( \Phi \circ \beta \right) \circ \left( \alpha \circ \Psi \right)(h') ) \leq \lambda + C$.
We then find a cycle of edges in $\Gamma^{(1)}$ from $h$ to $h'$ to $\left( \Phi \circ \beta \right) \circ \left( \alpha \circ \Psi \right)(h')$ to $\left( \Phi \circ \beta \right) \circ \left( \alpha \circ \Psi \right)(h)$ to $h$ of length
no more than $1 + 2K + \lambda + C$.  Let $M_2 = \max \left\{ L_2\left( \lambda + C \right), 1 + 2K + \lambda + C \right\}$.
A $2$-cell is glued onto $\Gamma^{(1)}$ along each $1$-cycle of length no more than $M_2$.  As the $H$ action on 
$\Gamma^{(0)}$ is transitive, we see that there are finitely many $H$ orbits of such cells.  This is $\Gamma^{(2)}$.

Let $u$ be a simple $1$-cycle in $\Gamma^{(1)}$ with vertices $x_1$, $x_2$, $\ldots$, $x_n$, $x_{n+1}=x_1$, and let
$v_i = (\alpha \circ \Psi)(x_i)$.  $d_{X^{(1)}}(v_i, v_{i+1}) \leq \lambda + C$, so there is a geodesic in
$X^{(1)}$ connecting $v_i$ to $v_{i+1}$, with length no more than $\lambda + C$.  Concatenating these paths results
in a $1$-cycle in $X$ with length bounded by $n\left(\lambda+C\right)$.  There is a filling in $X^{(2)}$ of this
$1$-cycle by $2$-cells.  The image under $\Phi \circ \beta$ of each of these $2$-cells results in a $2$-cell in
$\Gamma^{(2)}$ constructed as above.  These $2$-cells do not fill $u$, rather they form a filling of
`$\left( \Phi \circ \beta \right) \circ \left( \alpha \circ \Psi \right)(u)$'.  Consider the edge $[x_i, x_{i+1}]$
of $u$.  As above, the cycle from $x_i$ to $x_{i+1}$ to $\left( \Phi \circ \beta \right) \circ \left( \alpha \circ \Psi \right)(x_i)$ to $\left( \Phi \circ \beta \right) \circ \left( \alpha \circ \Psi \right)(x_{i+1})$ and back to $x_i$
has length no more than $1 + 2K + \lambda + C$, so it is a $2$-cell of $\Gamma^{(2)}$.  When these are combined with the
above $2$-cells lifted from $X^{(2)}$ a filling of $u$ is obtained.  It follows that $\Gamma{(2)}$ is simply connected.

Let $d^1_{\Gamma}$ be the Dehn function for filling $1$-cycles in $\Gamma$.
There is $L_3$ such that the boundary of each $3$-cell of $X^{(3)}$ contains no more than $L_3$ $2$-cells.  As above,
the image of each of these $2$-cells under $\Phi \circ \beta$ is itself a $2$-cell in $\Gamma^{(2)}$.  It follows
that given a $2$-cycle in $X^{(2)}$ with $n$ cells, the image under $\Phi \circ \beta$ is a 
$2$-cycle in $\Gamma^{(2)}$
with $n$ cells.  Let $\sigma$ be a $2$-cell of $\Gamma^{(2)}$, with vertices 
$x_1$, $x_2$, $\ldots$, $x_n$, $x_{n+1}=x_1$, and
let $y_i = \left( \Phi \circ \beta \right) \circ \left( \alpha \circ \Psi \right)(x_i)$.  $y_i$ can be joined to
$y_{i+1}$ by a geodesic of length no more than $\lambda + C$, so there is a cycle from $y_1$ to $y_2$ to $\ldots$
to $y_n$ to $y_1$ with length at most $M_2 \left( \lambda + C \right)$.  This $1$-cycle can be filled by fewer than
$d^1_{\Gamma}\left( M_2 \left( \lambda + C \right) \right)$ $2$-cells.  Denote this filling by $\sigma_b$.
Similarly there is a $1$-cycle obtained by traveling from $x_i$ to $x_{i+1}$ to $y_{i+1}$ to $y_i$ to $x_i$, with
length bounded by $1 + 2K + \lambda + C \leq M_2$, thus it can be filled by a single $2$-cell $\sigma_i$.
The $2$-cycle with top $\sigma$, bottom $\sigma_b$, and sides $\sigma_i$ has area no more than 
$M_2 + 1 + d^1_{\Gamma}\left( M_2 \left( \lambda + C \right) \right)$.  Let 
$M_3 = \max \{ L_3, M_2 + 1 + d^1_{\Gamma}\left( M_2 \left( \lambda + C \right) \right) \}$.  We glue a $3$-cell onto
a $2$-cycle $\omega$ in $\Gamma^{(2)}$ whenever $\omega$ has area no more than $M_3$.  This gives $\Gamma^{(3)}$.

Let $u$ be a connected $2$-cycle in $\Gamma^{(2)}$ composed of the $2$-cells $\nu_1$, $\nu_2$, $\ldots$,
$\nu_n$, $\nu_{n+1} = \nu_1$.  For each edge $[x,y]$ in the boundary of $\nu_i$, there
is a path $[ \left( \alpha \circ \Psi\right)(x), \left( \alpha \circ \Psi\right)(y)]$ in $X^{(1)}$ with
length no more than $\lambda + C$, and there are at most $M_2$ such edges for $\sigma_i$.  These paths
join to yield a $1$-cycle in $X$, which can be filled with $2$-cells of $X^{(2)}$.  Denote this filling
by $\left( \alpha \circ \Psi \right)( \nu_i )$.  The collection of each of these 
$\left( \alpha \circ \Psi \right)( \nu_i )$ spans a $2$-cycle in $X^{(2)}$, so it can be filled
by $3$-cells of $X^{(3)}$.  By construction, the image under $\Phi \circ \beta$ of each of these
$3$-cells constitutes a $3$-cell in $\Gamma^{(3)}$.  As above, each $\nu_i$ generates a $2$-cycle
with area bounded by $M_2 + 1 + d^1_{\Gamma}\left( M_2 \left( \lambda + C \right) \right)$ with
top $\nu_i$, bottom $(\nu_i)_b$, and the appropriate edges.  These also correspond to a single
$3$-cell each.  The $3$-complex spanned by the image under $\Phi \circ \beta$ of the fillings of the
$\left( \alpha \circ \Psi \right)( \nu_i)$, and the $3$-cells generated by each of the $\nu_i$ yield
a filling of $u$.  It follows that $\pi_2\left( \Gamma \right)$ is trivial. 

Suppose that we have constructed $\Gamma^{(k)}$ in such a way that for each $m \leq k$, there are finitely
many orbits of $m$ cells under the $H$ action, for $m < k$ $\pi_m\left( \Gamma \right)$ is trivial, and
there are constants $M_m$ such that the boundary of each $m$-cell in $\Gamma$ consists of no more than
$M_m$ $(m-1)$-cells.
Moreover, for $1 < m \leq k-1$ suppose that for an $m$-cell $\sigma$ in $X^{(m)}$, there is an $m$-cell
$\left( \Phi \circ \beta \right)(\sigma)$, which has vertices corresponding to the image of the 
vertices of $\sigma$ under $\left( \Phi \circ \beta \right)$ with geodesic edges as constructed above
for $2$ and $3$-cells, and that there is a constant $N_m$ such that for an $m$-cell $\sigma$, there is an
$m$-complex $\left( \Phi \circ \beta \right) \circ \left( \alpha \circ \Psi \right)(\sigma)$ having as vertices
the image under $\left( \Phi \circ \beta \right) \circ \left( \alpha \circ \Psi \right)$ of the vertices of $\sigma$
as the $2$-complex $\sigma_b$ was constructed above.

Let $u$ be an $k$-cell in $X^{(k)}$.  Denote the $(k-1)$-cells in the boundary of $u$ by 
$\sigma_1$, $\sigma_2$, $\ldots$, $\sigma_n$, $\sigma_{n+1} = \sigma_1$, with $n \leq L_k$.
Each $\sigma_i$ maps to $\left( \Phi \circ \beta \right)( \sigma_i)$, a $(k-1)$-cell in $\Gamma^{(k-1)}$, 
yielding a $(k-1)$-cycle in $\Gamma$ of area no more than $L_k$.

Let $\sigma$ be a $(k-1)$-cell of $\Gamma^{(k-1)}$.  The $(k-1)$-complex 
$\sigma_b = \left( \Phi \circ \beta \right) \circ \left( \alpha \circ \Psi \right)(\sigma)$ 
has area no more than $N_{k-1}$. For each $(k-2)$-cell, $\nu$, in the boundary of $\sigma$, there corresponds 
a $(k-2)$-subcomplex $\left( \Phi \circ \beta \right) \circ \left( \alpha \circ \Psi \right)(\nu)$ of $\sigma_b$,
which generate the boundary of $\sigma_b$.  We refer to $\nu$ and the
$(k-2)$-subcomplex $\left( \Phi \circ \beta \right) \circ \left( \alpha \circ \Psi \right)(\nu)$ as being parallel
cells.  Each vertex of $\nu$ is within an $H$-distance of $K$ from the corresponding vertex of 
$\left( \Phi \circ \beta \right) \circ \left( \alpha \circ \Psi \right)(\nu)$.  The $(k-2)$-subcomplex with top
$\nu$, bottom $\left( \Phi \circ \beta \right) \circ \left( \alpha \circ \Psi \right)(\nu)$, and sides obtained by
connecting the corresponding vertices, bounds a single $(k-1)$-cell by construction.  The $(k-1)$-complex
with top $\sigma$, bottom $\sigma_b$, and sides the $(k-1)$-cells so obtained, consists of at most
$1 + N_{k-1} + M_{k-1}$ $(k-1)$-cells.  Let $M_k = \max \{ L_k, 1 + N_{k-1} + M_{k-1} \}$.  We glue
a $k$-cell onto each $(k-1)$-cycle of $\Gamma^{(k-1)}$ with at most $M_k$ $(k-1)$-cells.  This yields
$\Gamma^{(k)}$.

Let $u$ be a connected $(k-1)$-cycle in $\Gamma^{(k-1)}$ composed of the $(k-1)$-cells $\nu_1$, $\nu_2$, $\ldots$,
$\nu_n$.  For each $\nu_i$ we form $\left( \alpha \circ \Psi \right)(\nu_i)$ as above.
The collection of all $\left( \alpha \circ \Psi \right)(\nu_i)$ forms a $(k-1)$-cycle,
denoted by $\left( \alpha \circ \Psi \right)(u)$ in $X^{(k-1)}$.  There is a filling in $X^{(k)}$
of $\left( \alpha \circ \Psi \right)(u)$ by $k$-cells $\eta_1$, $\eta_2$, $\ldots$, $\eta_j$.
Each $\left( \Phi \circ \beta \right)( \eta_i )$ gives a $k$-cell in $\Gamma$.  For each $\nu_i$, the
$(k-1)$-cycle with top $\nu_i$, bottom $(\nu_i)_b$, and the appropriate sides obtained by connecting
the parallel faces, bounds a $k$-cell, $\omega_i$.  The $k$-subcomplex spanned by the
$\left( \Phi \circ \beta \right)( \eta_i )$ and the $\omega_i$ forms a filling of $u$ in $\Gamma^{(k)}$.
As such, $\pi_k\left( \Gamma \right)$ is trivial.

The result follows by induction.
\endproof

\begin{cor}
A group quasiisometric to an $HF^\infty$ group with the isocohomological property, itself is
an $HF^\infty$ group with the isocohomological property.  In particular, in the class of $HF^\infty$
groups, the isocohomological property is a coarse invariant.
\end{cor}
\proof
By \cite{alonso} quasiisometric groups have equivalent higher Dehn functions.  The result follows from
Theorem \ref{thmClassification} and Lemma \ref{lemHFinfty}.
\endproof

As all polynomial growth groups are virtually nilpotent we obtain the following corollary.
\begin{cor}
The higher Dehn Functions for groups of polynomial growth are polynomially bounded.
\end{cor}

\section{Combing relatively hyperbolic groups}\label{sectCombRHG}
The concept of a relatively hyperbolic group was proposed by Gromov \cite{Gr2}, as a generalization
of hyperbolicity, in terms of a group action on a $\delta$-hyperbolic metric space.  Farb \cite{farb}
defined a more general `weak relative hyperbolicity' in terms of the geometry of a coned-off Cayley
graph, and introduced the bounded coset penetration property to allow estimates.  Bowditch \cite{bowditch}
defined a notion of relative hyperbolicity in terms of a group action on a graph and showed that his version
is equivalent to both Gromov's relative hyperbolicity and Farb's weak relative hyperbolicity with the 
bounded coset penetration property.  There are also notions of relative hyperbolicity due to Osin
\cite{Os2} and Mineyev-Yamen \cite{MY}, which are given in terms of relative Dehn functions and relative
complexes respectively, and which are equivalent to those above.  We will use the characterization of relative
hyperbolicity given by Farb, as weak relative hyperbolicity with the bounded coset penetration property.

Let $G$ be a finitely generated group endowed with the word-metric $d_G$.
A path in $G$ is an eventually constant function $p : \N \rightarrow G$ such that for each $n$,
$d_G( p(n), p(n+1) ) \leq 1$.  By identifying $G$ with the vertices of its Cayley graph,
we consider the path as a continuous function from $\R_+$ to $G$, with the interval $(m,m+1)$ mapped
onto an edge if $p(m) \neq p(m+1)$, or mapped to the vertex $p(m) = p(m+1)$.
The length of the path $p$ is $\sum_{n \in \N} d_G( p(n), p(n+1) )$.
As the path is eventually constant, it has finite length.  We say that $p$ starts at $p(0)$ and ends
at $\lim_{t \to \infty} p( t )$.

\begin{defn}
A combing of $G$ is a family $\sigma$ of paths, one for each element of $G$, satisfying the following.
\begin{enumerate}
\item[(1)] For each $g \in G$, $\sigma_g$ is a path starting at $e$ and ending at $g$. 
\item[(2)] There is a $K > 0$ such that for all $g$ and $g' \in G$ and all $t \in \N$, 
			$d_G( \sigma_g(t), \sigma_{g'}(t) ) \leq K d_G(g,g')$.
\end{enumerate}
$\sigma$ is a polynomial combing if there is a polynomial $P$ such that for all $g \in G$
the length of $\sigma(g)$ is no more than $P( d_G(g,e) )$.
\end{defn}

Of course this definition makes sense for any discrete metric space, once a distinguished basepoint has
been identified.  All automatic, hyperbolic, semihyperbolic, and CAT(0) groups are combable, with quasigeodesic
paths.

\begin{defn}
A combing $\sigma$ of a discrete metric space $(X,d)$ is said to have uniformly bounded return if there is a constant $N > 0$
such that for any $x \in X$ and any $y \neq x \in X$, the set $\{ t \, | \, \sigma_x(t) = y \}$ has cardinality bounded by $N$.
\end{defn}
In particular any geodesic or quasigeodesic combing has uniformly bounded return.  Our main interest in this property
is the following.

\begin{lem}\label{timebdd}
Let $(X,d)$ be a discrete metric space with bounded geometry, endowed with a combing of uniformly bounded return $\sigma$.
There is a function $f$ with the following property:
For any $x$,$y \in X$, with $t_x$ and $t_y$ the times that $\sigma_x$ and $\sigma_y$ settle
at $x$ and $y$ respectively, then $|t_y - t_x| \leq f( d(x,y) )$.
\end{lem}
\proof
Suppose that $K$ is the combing constant of $\sigma$, and let $N$ be the return constant.  We suppose without loss of
generality that $t_y \leq t_x$.  For all $t_y \leq t \leq t_x$, $d( y, \sigma_x(t) ) \leq K d(x,y)$ so that $\sigma_x(t)$
is in the ball of radius $K d(x,y)$ centered at $y$.  Let $V( r )$ be the volume of a ball of radius $r$.  Then there are
at most $V( K d(x,y) )$ possibilities for $\sigma_x(t)$.  As $\sigma_x$ can take each of these values at most $N$ times,
it must settle at $x$ by $t_y + N V( K d(x,y) )$.  Take $f(r) = N V( K r )$.
\endproof

Let $\Gamma$ be the Cayley graph of the group $G$ with respect to some finite generating set,
and let $H_1$, $H_2$, $\ldots$, $H_n$ be subgroups.  
The coned-off Cayley graph $\hat G$ is obtained from $\Gamma$ by adding one vertex for each coset 
$gH_i$ of each $H_i$ in $G$, and adding
an edge of length $1/2$ between each element vertex and the vertex of each coset to which the element belongs.

$\hat G$ is said to satisfy the bounded coset penetration (BCP) property if for each
constant $k$ there exists a constant $c = c(k)$ such that for every
pair of $k$-quasigeodesics $p$ and $q$ in the coned-off Cayley graph with 
the same endpoints and without backtracking, satisfy

\begin{enumerate}
\item[(1)] If $p$ penetrates a coset $gH_i$ and $q$ does not penetrate 
		$gH_i$, then the point at which $p$ enters the coset is
		at most a $\Gamma$-distance of $c(k)$ from the point at which $p$ leaves
		the coset.

\item[(2)] If $p$ and $q$ both penetrate a coset $gH_i$, then the points
		at which $p$ and $q$ enter $gH_i$ are at most a $\Gamma$-distance
		of $c(k)$ from each other.  Similarly the points at which $p$ and $q$ exit
		the coset are within a $\Gamma$-distance of $c(k)$ from one-another.  
\end{enumerate}

\begin{defn}
$G$ is relatively hyperbolic with respect to the subgroups $H_1$, $H_2$, $\ldots$, $H_n$, if 
$\hat G$ is $\delta$-hyperbolic for some $\delta > 0$, and satisfies the BCP property.
\end{defn}

\begin{lem}
Let $\hat G$ be the coned-off Cayley graph of the group $G$ with respect to the family of subgroups $H_1$, $H_2$, $\ldots$,
$H_n$.  There is a family of geodesics $\alpha_g$ in $\hat G$,
for $g \in G$, with $\alpha_g(0) = e$ and $\alpha_g(\ell_{\hat G}(g))=g$, such that if
$\alpha_g(t') = \alpha_{g'}(t')$ for some $t'$, then $\alpha_g(t) = \alpha_{g'}(t)$ for all $t \leq t'$.
\end{lem}
\proof
Enumerate the group elements as $e = g_0$, $g_1$, $g_2$, \ldots, with 
$\ell_G(g_i) \leq \ell_G(g_{i+1})$ for all $i$.  Let $[e,g_i]$ be some fixed geodesic in $\hat G$ from $e$ to $g_i$.
As geodesics, for each positive integer $t$, $[e, g_i](t)$ is a group element vertex, not a coset vertex.
Clearly $[e,g_0]$ is the trivial path.
If $\ell_G(g_i) = 1$ then $[e,g_i]$ is the geodesic consisting of the one edge from $e$ to $g_i$.
If $\ell_{\hat G}(g_i) = 1$ with $\ell_G(g_i) \neq 1$ then $g_i$ is in the identity coset, but does not differ from
$e$ by a generator, so the geodesic from $e$ to $g_i$ must penetrate the coset $H$ and then travel to $g_i$.
In each of these cases we let $\alpha_{g_i} = [e, g_i]$.  We inductively construct the remainder as follows.

Suppose we have constructed $\alpha_{g_{i'}}$ for all $i' < i$.  Define 
\[l_i(j) = \max \left\{ t \, | \, [e,g_i](t) = \alpha_{g_j}(t) \right\}\] 
for $j = 0, 1, \ldots, i-1$.  Let $j$ be the largest of these indices which maximize $l_i$, with $l_i(j) = \bar t$.  
We define $\alpha_{g_i}(t) = \alpha_{g_j}(t)$ for $t \leq \bar t$ and $\alpha_{g_i}(t) = [e, g_i](t)$
for $t > \bar t$.

The final claim is obvious for $\ell_{\hat G}(g) \leq 1$ or $\ell_{\hat G}(g') \leq 1$.
Suppose that $\alpha_{g_i}(t_0) = \alpha_{g_i'}(t_0)$ for some $t_0 > 0$, with $i > i'$.  
We take $t_0$ to be the largest time at which these two paths intersect.  Take $j$ and $\bar t$ 
as in the construction of $\alpha_{g_i}$ above.  If $j = i'$ then $t_0 = \bar t$ and the two paths
agree for all time less than $t_0$.  Otherwise, $\bar t \geq t_0$, so for all time less than 
$t_0$ $\alpha_{g_i}(t) = \alpha_{g_j}(t)$.  Moreover, $\alpha_{g_j}$ intersects $\alpha_{g_{i'}}$ at $t_0$
with both $j$ and $i'$ strictly less than $i$.  By induction they agree for all time less than $t_0$.
Therefore $\alpha_{g_i}$ agrees with $\alpha_{g_{i'}}$ as claimed.
\endproof

As $\hat G$ is hyperbolic, this set of paths forms a geodesic combing.  
Suppose that each $H_i$ is combable with uniformly bounded return, with length bounded by polynomial $P$, and
denote the combing of $H_i$ by $\sigma^i$.

The combing of the subgroups and the combing $\alpha$ of $\hat G$ induce a system of paths $\beta_g$ in $G$ as follows.
For $g \in G$ consider the path $\alpha_g$ in $\hat G$.
In each unit interval of time $[t, t+1]$, for $t \in \N$, where $\alpha_g$ is not stationary, $\alpha_g$ either travels
along one edge of length $1$, or two edges of length $1/2$.  In the latter case, the $\alpha_g$ also penetrates a coset
during this time interval.
If during this interval $\alpha_g$ travels from $a$ to $b$ along a single edge of length $1$, then $d_G( a, b) = 1$
so $a$ and $b$ are joined by an edge in $\Gamma$.  The corresponding path $\beta_g$ will travel in $\Gamma$ from 
$a$ to $b$ along this edge, and stay at $b$ for $P( c(1) ) - 1$ additional units of time.  On the other hand,
if $\alpha_g$ penetrates the coset vertex $aH_i$ to reach $b$, then $b = ah$ for some $h \in H_i$.  The
corresponding path $\beta_g$ will travel in $\Gamma$ from $a$ to $b$ along the path $a \sigma^i_h$, the path
$\sigma^i_h$ translated to start at $a$.  If this portion of $\beta_g$ takes less than $P( c(1) )$ units of
time to traverse, then the path remains stationary at $b$, until the full $P( c(1) )$ time has expired, while 
if this portion takes longer than $P(c(1))$ to traverse, it is not altered.
$\beta_g$ is the concatenation of these paths.
In this way the projection of $\beta_g$ into $\hat G$ is exactly $\alpha_g$.

\begin{thm}\label{thmPolyComb}
Suppose that $G$ is relatively hyperbolic with respect to the subgroups $H_i$.
The system of paths $\beta_g$ constructed above is a polynomial length combing of $G$. 
\end{thm}

It is not surprising that a group relatively hyperbolic to combable subgroups is itself combable; however,
that the combing has polynomially bounded length when the subgroup combings do is new, and relies on the recent
result that the bounded coset penetration function $c$ is itself polynomially bounded \cite{jor}.
We are curious about how to construct a combing without using estimates based on $c$.  In particular,
if the subgroups are quasigeodesically combable, must the full group be quasigeodesically combable?

The proof of Theorem \ref{thmPolyComb} consists of a series of lemmas.  We fix two elements $g$ and $g'$ of $G$, with
$d_G( g, g') = 1$.  $\alpha_g$ and $\alpha_{g'}$ are two $1$-quasigeodesics which start at the same group
element and end at group elements at distance $1$ apart.  We will use bounded coset penetration
to compare them.  
\begin{defn}
Two cosets $aH_i$ and $bH_j$ are said to be synchronous with respect to $g$ and $g'$
if for some positive integer $t$, $\alpha_g(t)$ and $\alpha_{g'}(t)$ are elements in $aH_i$ and $bH_j$
respectively.
\end{defn}  

\begin{lem}\label{lem3.7}
Let $aH_i$ and $bH_j$ be synchronous cosets for $g$ and $g'$.  Assume $\alpha_g$ enters $aH_i$ at $a$ and exits
at $ah$, while $\alpha_{g'}$ enters $bH_j$ at $b$ and exits at $bh'$.  There is a constant $M$, independent
of the cosets and the group elements, for which $d_G(a,b) \leq M$ and $d_G( ah, bh') \leq M$.
\end{lem}
\proof
If $aH_i = bH_j$ then $a=b$, $d_G( a, ah ) \leq c(1)+1$, and $d_G( b, bh') \leq c(1)+1$.  
Otherwise, $aH_i \neq bH_j$, 
so $\alpha_{g'}$ can not intersect the coset $aH_i$, and $\alpha_g$ cannot intersect $bH_j$.
Therefore, $d_G( a, ah ) \leq c(1) + 1$ and $d_G( b, bh' ) \leq c(1) + 1$.  It is sufficient
to show the result only for the exiting points.

If $K$ is the combing constant for $\alpha$, then $d_{\hat G}( ah, bh' ) \leq K+1$.
Let $u$ be a geodesic in $\hat G$ connecting $ah$ to $bh'$, let $\gamma$ be the path
from $ah$ to $g$ followed by $\alpha_g$, and let $\gamma'$ the path from $bh'$ to $g'$ followed
by $\alpha_{g'}$.  Denote the path from $g'$ to $g$ obtained by following the edge connecting
the two points as $s$.  

Suppose that $u$ does not penetrate any of the same cosets at $\gamma$ or $\gamma'$.
Then we have a geodesic $\gamma$ from $ah$ to $g$ and a $(2K+3)$-quasigeodesic, obtained by
concatenating $u$, $\gamma'$, and  $s$, also from $ah$ to $g$, without backtracking, both
starting and ending at the same point.
By BCP any coset penetrated by $u$ can be traveled through a distance at most $c( 2K+3 )$.
It follows that $d_G( ah, bh' ) \leq (K+1) c( 2K+3 )$.

Suppose that $u$ penetrates one of the cosets that $\gamma'$ penetrates, but not none that are 
penetrated by $\gamma$.  We assume that $u$ is picked
so that $u$ and $\gamma'$ agree after they meet.  Let $u$ be the concatenation of $u'$ and $u''$, where $u'$
is a geodesic from $ah$ to the point where $u$ and $\gamma'$ first meet, and $u'$ is the portion of $u$ where
it follows $\gamma'$.  Let $\omega$ be the portion of $\alpha_g$ from $e$ to $ah$, and $\omega'$ the portion of
$\alpha_{g'}$ from $e$ to the point where $\gamma'$ and $u$ first meet.  $\omega'$ is a geodesic and the 
concatenation of $\omega$ and $u'$ is a $(2K+2)$-quasigeodesic without backtracking, both starting and ending at
the same point.  By BCP $u'$ can not penetrate any coset by more than $c( 2K+2 )$, so the $G$-distance
between $ah$ and the endpoint of $u'$ is no more than $(K+1)c(2K+2)$.  As $\alpha_{g'}$ and $u$ agree past this
endpoint, $u''$ can penetrate any coset by at most $c(1)$.  The $G$-distance between $bh'$ and the endpoint
of $u'$ is thus bounded by $(K+1) c(1)$.  It follows that the $G$-distance from $ah$ to $bh'$ is no more than
$(K+1) ( c(2K+2) + c(1) )$.  The case where $u$ penetrates some of the same cosets as $\gamma$, but none of those 
penetrated by $\gamma'$ is similar.

Suppose that $u$ penetrates cosets penetrated by both $\gamma$ and $\gamma'$.
We assume that $u$ is the concatenation of $u'$, $u''$, and $u'''$, where $u'$ is the portion of $\gamma$ between
$ah$ and the last point where $u$ meets $\gamma$, and $u'''$ is the portion of $\gamma'$ between the first point
where $u$ meets $\gamma'$ and $bh'$.  In this case, $u''$ penetrates no coset which is penetrated by either
$\gamma$ or $\gamma'$.  By BCP as above, the $G$-distance between the endpoints of $u''$ is no more than
$(K+1) c( 2K+3 )$.  The $G$-distance between the endpoints of $u'$ is no more than
$(K+1) c( 1 )$, and the $G$-distance between the endpoints of $u'''$ is no more than $(K+1)c(1)$.
The $G$-distance from $ah$ to $bh'$ is therefore no more than $(K+1)( c( 2K+3 ) + 2 c(1) )$.
\endproof

\begin{lem}\label{lem3.8}
Let $aH_i$ and $bH_j$ be synchronous cosets for $g$ and $g'$.  Assume $\beta_g$ enters $aH_i$ at time $t_g$ and exits
at $\bar t_g$, while $\beta_{g'}$ enters $bH_j$ at time $t_{g'}$ and exits at $\bar t_{g'}$.  There is a constant $T$,
independent of the cosets and the group elements, for which 
$| t_g - t_{g'} | \leq T$ and $| \hat t_g - \hat t_{g'} | \leq T$.
\end{lem}

\proof
Let $x = \alpha_g( \hat t ) = \alpha_{g'}( \hat t )$ be the last group element where the two paths $\alpha_g$ and
$\alpha_{g'}$ meet, and let $t_x$ be the time for which $\beta_g(t_x) = \beta_{g'}(t_x) = x$.  
For all times $t \leq t_x$, $\beta_g(t) = \beta_{g'}(t)$ so the two paths enter and exit all earlier cosets at exactly the same time.

If both $\alpha_g$ and $\alpha_{g'}$ penetrate the coset $xH$, then they exit the coset $xH$ at points with $G$-distance
no more than $c(1)$ apart, while they enter $xH$ at the same time.  As the combing of $H$ has uniformly bounded return,
the times at which $\beta_g$ and $\beta_{g'}$ leave $xH$ differ by no more than $f( c(1) )$, where $f$ is the function
from Lemma \ref{timebdd}.   
Otherwise, neither $\alpha_g$ nor $\alpha_{g'}$ travel in $xH$ by more than $c(1)$, so the time $\beta_g$ and $\beta_{g'}$ spend
in $xH$ is exactly $P(c(1))$, so they exit at exactly the same time.

If $aH_i \neq bH_j$, then those cosets must be reached after time $t_x$.  
As $\alpha_{g'}$ does not penetrate $aH_i$ and $\alpha_{g}$ 
does not penetrate $bH_j$, then $\beta_g$ and $\beta_{g'}$ spend exactly $P(c(1))$ units of time in $aH_i$ and $bH_j$ respectively.
Thus $\beta_g$ and $\beta_{g'}$ exit cosets $aH_i$ and $bH_j$ with the same time difference with which they were entered.
\endproof

The previous lemmas show that the paths $\beta$ form a synchronous combing.  It remains to show that they are of 
polynomial length.

\begin{lem}
The paths $\beta$ are of polynomial length.
\end{lem}
\proof
Let $g \in G$.  Let $\gamma$ be a geodesic in $G$ from $e$ to $g$, and let $\hat \gamma$ be the projection into $\hat G$.
$\hat \gamma$ is a $\ell_G( g ) + 1$-quasigeodesic, starting and stopping at the same group elements.
If $aH$ is a coset penetrated by $\alpha_g$ but not by $\hat \gamma$, then $\alpha_g$ can travel no more than
$c( \ell_G(g)+1 )$ through $aH$.  Then $\beta_g$ can travel no more than $P( c( \ell_G(g) + 1 )$ in the coset $aH$.
If $\alpha_g$ and $\hat \gamma$ both penetrate $aH$, then they enter and exit
within $c( \ell_G(g)+1 )$ of each other.  As $\hat \gamma$ travels no more than $\ell_G(g)$ through $aH$,
then $\alpha_g$ can travel no more than $\ell_G(g) + 2 c( \ell_G(g) + 1 )$ through the coset.  Then $\beta_g$
can travel no more than $P( \ell_G(g) + 2 c( \ell_G(g) + 1 ) )$ through $aH$.  As $\alpha_g$
travels through no more than $\ell_G(g)$ cosets, the length of $\beta_g$ is no more than
$\ell_G(g)  P( \ell_G(g) + 2 c( \ell_G(g) + 1 ) )$.  By Proposition 2.2.7 of \cite{jor}, the coset penetration function 
$c$ is polynomially bounded.  
\endproof

%

\section{Classifying spaces for relatively hyperbolic groups}
The characterization given in Theorem \ref{thmClassification} raises natural questions for a group $G$ 
relatively hyperbolic to a family of subgroups $H_i$.  If each of the $H_i$ are $HF^\infty$ groups, is $G$?  
And if so, if the classifying spaces for $H_i$ all have polynomially 
bounded higher Dehn functions, does the one for $G$?  
As polynomial combability implies both of these properties, the results of the previous section suggest that 
these are possible.  The main result of this section is the following.
\begin{thm}
Let $G$ be relatively hyperbolic to $HF^\infty$ subgroups $H_1$, $H_2$, $\ldots$, $H_n$.
If each $H_i$ is isocohomological, then $G$ is isocohomological.
\end{thm}
This will be accomplished by showing that if each $H_i$ is $HF^\infty$ then so is $G$, and that
if each $H_i$ has the appropriate Dehn functions, so does $G$.  We do this by constructing out of the
appropriate spaces for the $H_i$, a classifying space for $G$, similar to Gersten's construction in \cite{ger}.  
We remark that in \cite{dah}, there is a construction of a classifying space for $G$ when the subgroups 
are torsion-free.  In our construction this requirement is unnecessary.  For notational convience we will
work with only a single subgroup.  The proof for many subgroups is virtually identical.

Assume that $G$ is relatively hyperbolic with respect to $H$, $H$ has, as the universal cover of 
its classifying space, $X$, an aspherical polyhedral 
complex with finitely many orbits of cells in each dimension, and has weighted Dehn functions 
$d_{w,H}^n$.  We may assume that $X$ contains a copy of the Cayley graph of $H$.

Let $\Gamma$ be the polyhedral complex constructed as follows.
The vertices of $\Gamma$ are precisely the vertices of $X$ and its translates by elements of $G$.

The edges of $\Gamma$ are the edges of $X$ and those of its translates.  There is also an edge between
two group elements $g$ and $g'$ if $d_G(g,g') = 1$.  This way we are assured that a copy of the Cayley graph
of $G$ lies in $\Gamma$.

Fix a system of paths $\sigma$ in $H$, such that $\sigma_h$ is a geodesic in $H$ connecting the identity of $H$ to $h$.
This system need not be a combing, but the systems $\alpha$ and $\beta$ can be constructed using $\sigma$ as 
in section \ref{sectCombRHG}. 
$\alpha$ is still a combing of the coned-off graph, but $\beta$ need not be a combing.
The results of Lemma \ref{lem3.7} and of Lemma \ref{lem3.8} are still valid, as they did not rely on $\sigma$ being 
a combing, but only on the structure of $\alpha$ and the bounded coset penetration property.
Thus if $g$ and $g'$ are two elements of $G$, with $d_G( g, g' ) = 1$, we consider the paths 
$\beta_g$ and $\beta_{g'}$.  If $g$ and $g'$ are not in the same coset, 
For any $t$, $\beta_g(t)$ and $\beta_{g'}(t)$ are in synchronous cosets, say $aH$ and $bH$.  Assume that $\beta_g$ and
$\beta_{g'}$ enter cosets $aH$ and $bH$ at vertices $a$ and $b$ at times $t_a$ and $t_b$ respectively.  
By Lemma \ref{lem3.7}, $d_G( a, b ) \leq M$, while by Lemma \ref{lem3.8} $| t_a - t_b | \leq T$.  
As $\beta_g$ and $\beta_{g'}$ travel no more than $c(1)$ through $aH$ and $bH$ 
(unless $aH = bH$, but then $a = b$ by construction ) $d_G( \beta_g(t), \beta_{g'}(t) ) \leq M + T + 2 c(1)$.

Otherwise, $g$ and $g'$ are in the same coset.  If $\beta_g$ and $\beta_{g'}$ enter the coset at the same point, 
$\beta_g(t) = \beta_{g'}(s)$, then by construction $t = s$ and for all $t' < t$, $\beta_g(t') = \beta_{g'}(t')$.
In case $\beta_g$ and $\beta_{g'}$ do not enter the coset at the same point.  Assume that $\beta_g$ and $\beta_{g'}$
enter the coset at points $\bar g$ and $\bar{g'}$ at times $t_g$ and $t_{g'}$ respectively.  Assume also that
$t_g \leq t_{g'}$.  As above for all $t \leq t_{g'}$, $d_G( \beta_g(t), \beta_{g'}(t) ) \leq M + T + 2 c(1)$.

The construction given above allows us to `comb' an edge with vertices lying in different translates of $X$,
as a cone to the basepoint.
The $2$-cells of $\Gamma$ are the $2$-cells of $X$ and those of its translates, as well as some special $2$-cells
glued along a $1$-cycle consisting of no more than $2M+2T+4c(1)+2$ edges, all of whose vertices do not lie in the
same translate of $X$.  We call such $2$-cells the {\it{mixing}} $2$-cells, as they are the only $2$-cells with 
vertices from multiple translates of $X$.

It is known by the work of Farb in \cite{farb} that the usual Dehn function for $G$ is the maximum of the Dehn functions for the $H_i$.  Rather than appealing to Lemma \ref{lemBddWeightedDehn}, we estimate the weighted Dehn function directly
as follows. 
\begin{lem}
$\Gamma^{(2)}$ is simply connected and has weighted Dehn function for filling $1$-boundaries bounded by 
$f(x) = x(x+1)\left(d_{w,H}^1(x)+1\right)$.
\end{lem}
\proof
As $d_{w,H}^*$ is defined in terms of a length function on $H$; Different length functions yield different
weighted Dehn functions, with quasiisometric lengths yielding equivalent functions.  As weighted 
Dehn functions
are defined only up to this equivalence, we assume that  $d_{w,H}^1$ is a weighted Dehn function corresponding
to the length $\ell_{X}( v ) = d_{\Gamma}(e, v)$, where $e$ is the vertex of $\Gamma$ corresponding
to the identity element of $G$.  This length is quasiisometric to the length on the $1$-skeleton of $X$, so
$d_{w,H}^1$ is equivalent to a weighted Dehn function on $X$.
Let $u$ be a $1$-boundary in $\Gamma^{(1)}$.  Label the vertices of $u$ as $v_1, v_2, \ldots, v_n$.
Assume that $v_i$ does not correspond to a group element, lying in the translate $g X$.
From the finitely many orbits property of $X$, there is a constant $L$ such that every vertex 
in $\Gamma^{(0)}$ not corresponding to a group element,
is within a $H$-distance of $L$ from a group element vertex in the same translate of $X$.  
Let $v'_i$ be a group element vertex in the same translate as $v_i$, within an $H$-distance of $L$.
When $v_j$ does correspond to a group element, let $v'_j = v_j$.
If $v_i$ was a non-group element vertex, $v'_{i-1}$ and $v'_{i+1}$ are in the same translate as $v'_i$,
as only group element vertices have edges exiting their translate of $X$.  Thus there is an
edge path in $X$ connecting $v'_{i-1}$ to $v'_i$, with only group element vertices.  As 
$d_{gX}(v'_{i-1}, v'_i) \leq 2L+1$, and $X$ is quasiisometric to the Cayley graph of $H$, these paths can
be chosen to have length no more than $N$, for some universal constant $N$, independent of $u$.  Similarly
$v'_i$ can be connected to $v'_{i+1}$ through a path of edge of length at most $N$, passing through
only group element vertices.  Doing this for all non-group element vertices on $u$ results in a new
cycle of edges $u'$ with length no more than $N$ times that of $u$, which travels through only group
element vertices.  This transition from $u$ to $u'$ in $\Gamma$ is the result of passing portions of $u$
through $2$-cells of $gX$, so that $u$ and $u'$ are homotopic in $\Gamma^{(2)}$, thus if $u'$ is contractible
so is $u$.  Moreover $| u' |_w \leq (2L+1)^2 | u |_w$.

Thus assume that $u$ contains only vertices corresponding to group elements.  Again, denote these
by $v_1$, $v_2$, $\ldots$, $v_n$, $v_{n+1} = v_1$, and consider the $\beta_{v_i}$ paths.  If $v_i$ and $v_{i+1}$
lie in different cosets, $d_G( \beta_{v_i}(t), \beta_{v_{i+1}}(t) ) \leq M+T+2c(1)$ as above, so 
the cycle from $\beta_{v_i}(t), \beta_{v_i}(t+1), \beta_{v_{i+1}}(t+1), \beta_{v_{i+1}}(t)$ has length
no more than $2M + 2T + 4c(1) + 2$, and thus defines a $2$-cell, and the region in $u$ between $\beta_{v_i}$
and $\beta_{v_{i+1}}$ is covered by at most $c(1) ( \ell_{\Gamma}( v_i ) + 1 )$ of these cells.  
Then assume $v_i$ and $v_{i+1}$ lie in the same coset.  Denote by 
$g_i$ the point where $\beta_{v_i}$ enters the coset $v_i H$.  As $\beta_{v_i}$ and $\beta_{g_i}$ agree
up to $g_i$, $\beta_{v_i}$ is the concatenation of $\beta_{g_i}$ with $\gamma_i$, where $\gamma_i$ is the path
$\beta_{v_i}$ travels through $v_i H$.  For all $t$ $d_{\Gamma}( \beta_{g_i}(t), \beta_{g_{i+1}}(t) ) \leq M + T + 2c(1)$, so
the cycle from $\beta_{g_i}(t), \beta_{g_i}(t+1), \beta_{g_{i+1}}(t+1), \beta_{g_{i+1}}(t)$ defines a $2$-cell,
and there are no more than $c(1) ( \ell_{\Gamma}( v_i ) + 1 )$ of these cells.   

Consider the cycle $g_i, v_i, v_{i+1}, g_{i+1}$, following $u$ from $v_i$ to $v_{i+1}$, 
$\gamma_j$ from $v_j$ to $g_j$, and a path within the coset from $g_i$ to $g_{i+1}$, consisting only of group 
element vertices,  which can be chosen to have length no more than $2c(1)+1$.
Denote this cycle by $\eta$. The weighted length, in $G$, of this cycle is bounded by 
$\left( 2c(1) + 1 \right)| u |_w$.  Let $g$ be a vertex of $\eta$ with minimal length in $G$.
Let $\eta'$ be $\eta$ translated by $g^{-1}$.  That is, $\eta'$ is the cycle in $X$,
$g^{-1} g_i, g^{-1} v_i, g^{-1} v_{i+1}, g^{-1} g_{i+1}$, following the translated paths.  
\begin{eqnarray*}
| \eta' |_w & = & \sum_{\sigma' \in \eta'} | \sigma' |_w \\
 & = & \sum_{\sigma' \in \eta'} \sum_{v' \in \sigma'} \ell_{\Gamma}( v' )\\
 & \leq & \sum_{\sigma \in \eta} \sum_{v \in \sigma}\left( \ell_{\Gamma}(g) + \ell_{\Gamma}( v )\right)\\
 & \leq & \sum_{\sigma \in \eta} \sum_{v \in \sigma}\left( \ell_{\Gamma}(g) + \ell_{\Gamma}( v )\right)\\
 & \leq & \sum_{\sigma \in \eta} \sum_{v \in \sigma}\left( 2 \ell_{\Gamma}( v ) + L \right)\\
 & \leq & (L + 2) | \eta |_w
\end{eqnarray*}

There is a filling, $\omega'$, of $\eta'$ with $| \omega' |_w \leq d_{w,H}^1( \eta' )$.  
Let $\omega$ be $\omega'$ translated by $g$, a filling of $\eta$ in $g X$.  If $w$ is a vertex in $\omega'$, then 
$\ell_{\Gamma}( gw ) \leq \ell_{\Gamma}(g) + \ell_{\Gamma}(w)$.  If $w_1, w_2, \ldots, w_j$ are
the vertices of a cell $\nu'$ in $\omega'$ then 
\begin{eqnarray*}
| g \nu' |_w & \leq & \sum_k \ell_{\Gamma}(g w_k ) \\
 & \leq & \sum_k \left( \ell_{\Gamma}(g) + \ell_{\Gamma}( w_k )\right) \\
 & \leq & j \ell_{\Gamma}(g) + \sum_k \ell_{\Gamma}( w_k ) \\ 
 & = & j \ell_{\Gamma}(g) + | \nu' |_w\\
 & \leq & J | \eta |_w + | \nu' |_w
\end{eqnarray*}
Where $J$ is a constant such that no $2$-cell in $X$ has more than $J$ vertices.

\begin{eqnarray*}
| \omega |_w & \leq & \sum_{\nu' \in \omega'} | g \nu' |_w\\
 & \leq & \sum_{\nu' \in \omega'} \left( J | \eta |_w  + | \nu' |_w \right)\\
 & \leq & J | \eta |_w | \omega' |_w + | \omega' |_w\\
 & \leq & \left( J | \eta |_w + 1 \right) d_{w,H}^1( |\eta'|_w )\\
 & \leq & \left( J | \eta |_w + 1 \right) d_{w,H}^1\left( (L+2) |\eta|_w \right))
\end{eqnarray*}
$u$ has at most $| u |$ possible such $\eta$ cycles.  These can be filled with weighted length
at most $| u |_w \left( J |u|_w + 1 \right) d_{w,H}^1\left( (L+2) |u|_w \right)$.  The remainder
can be filled by cells with weighted length at most $c(1) |u|_w ( | u |_w + 1 )$, as above.
As such, 
\begin{eqnarray*}
d_{w,G}^1( |u|_w ) & \leq & | u |_w \left( J |u|_w + 1 \right) d_{w,H}^1\left( (L+2) |u|_w \right)\\
 & & + c(1) |u|_w ( | u |_w + 1 ).  
\end{eqnarray*}
\endproof

\begin{cor}
If $X$ has polynomially bounded weighted $1$-dimensional Dehn function, so does $\Gamma$.
\end{cor}

Let $\Delta$ be a mixing $2$-cell. We `comb' $\Delta$ to the basepoint by combing each of the edges as above, and
make the appropriate identifications as collections of $2$-cells.  Bounding the area of each of these sides, we
uniformly bound the area of each layer of the cone in terms of the $1$-dimensional Dehn functions of $\Gamma$
and $X$.  An issue arises when an edge of $\Delta$ is not a mixing edge, but this is overcome using BCP
and the fact that one of the faces of $\Delta$ is mixing.  The details are as follows.

Assume that the vertices of $\Delta$ are $v'_1$, $v'_2$, $\ldots$, $v'_k$, $v'_{k+1}=v'_1$.
For each $i$, if $v'_i$ is a group element vertex let $v_i = v'_i$.  Otherwise, there is a group element vertex $v_i$
within an $X$-distance $L$ from $v'_i$.  Connect $v_i$ to $v_{i+1}$ via a shortest edge path through group element
vertices.  This length is no more than $\lambda(2L+1)+C$.  Denote this cycle by $u$.  It has length at most 
$\left(\lambda(2L+1)+C\right)\left( 2M + 2T + 4 c(1) + 2 \right)$.  For each vertex $v$ of $u$, denote the point
where $\beta_v$ enters the coset $vH$ by $v^*$.  $\beta_v$ splits as the concatenation $\beta_v = \beta_{v^*} \gamma_v$
with $\gamma_v$ the path traveled by $\beta_v$ from $v^*$ to $v$.  For $v$ and $w$ consecutive vertices on $u$,
there is a path connecting $v^*$ and $w^*$ passing through only group element vertices, with length no more than $M$.  
The distance from $\beta_{v^*}(t)$ to $\beta_{w^*}(t)$ is no more than $M + T + 2c(1)$.  Let $u^*_t$ be the path
obtained by connecting $\beta_{v^*}(t)$ to $\beta_{w^*}(t)$ by a geodesic path, for each pair of consecutive
vertices $v$ and $w$ of $u$.  The distance from $\beta_{v^*}(t)$ to $\beta_{v^*}(t+1)$ is no more than $1$,
so consider the `drum' with top $u^*_{t+1}$, bottom $u^*_t$, and sides the cycles from $\beta_{v^*}(t)$ to 
$\beta_{v^*}(t+1)$ to $\beta_{w^*}(t+1)$ to $\beta_{w^*}(t)$.  The discussion above shows that each of these 
sides constitute a single $2$-cell, while $u^*_t$ and $u^*_{t+1}$ have lengths no more than 
$W = \left( M + T + 2c(1) \right)\left( \lambda(2L+1)+C \right)\left( 2M + 2T + 4 c(1) + 2 \right)$.
If $d^1$ denotes the Dehn function for filling $1$-dimensional cycles in $\Gamma$, then $u^*_t$ and $u^*_{t+1}$
can each be filled by at most $d^1(W)$ cells.  The `drum' has surface area no more than $W + 2 d^1(W)$.

If all vertices of $u$ lie in the same translate of $X$, say $gX$, then the portion of the cone above
$u^*$ with bottom face bounded by $u^*$, top face bounded by $u$, and sides bounded by the cycles $v$ to $w$
to $w^*$ to $v^*$ to $v$, for $v$ and $w$ consecutive vertices of $u$, is a $2$-boundary in $gX$, so can be 
filled using only those $3$-cells from the $gX$ structure.  Otherwise, not all vertices
of $u$ lie in $gX$. 
Let $w_1$, $w_2$, $\ldots$, $w_j$ be a maximal string of consecutive vertices of $u$ all lying in $gX$.  
There is a vertex $v$ of $u$ not lying in $gX$.  If $\beta_v$ does not penetrate $gH$, no $\beta_{w_i}$
can penetrate $gH$ by more than $c(1) \left( \lambda(2L+1)+C \right)$, so 
$d( w_i, w^*_i ) \leq c(1) \left(\lambda(2L+1)+C\right)\left( 2M + 2T + 4 c(1) + 2 \right)$.  Let $\gamma_i$ be 
a shortest length path in $gX$ from $w_i$ to $w^*_i$ through only group element vertices.
There is a path in $gX$ connecting $w_1$ to $w_j$ through group element vertices, with length
bounded by $\lambda \left(\lambda(2L+1)+C\right)\left( 2M + 2T + 4 c(1) + 2 \right) + C$.  
There is a similar path of length no more than $\lambda W + C$ connecting $(w_1)^*$ to $(w_j)^*$.
Consider the $2$-chain  with faces $v$ to $w_i$ to $w_{i+1}$, $v^*$ to $w^*_i$ to $w^*_{i+1}$, $v$ to $w_i$ to 
$w^*_i$ to $v^*$, $v$ to $w_{i+1}$ to $w^*_{i+1}$ to $v^*$, and $w_i$ to $w^*_i$ to $w^*_{i+1}$ to $w_{i+1}$.  As above these lengths are all bounded by a positive constant $K$, independent of $u$, so this wedge has surface area
no more than $5d^1(K)$.  
If $w$ and $w'$ are two adjacent vertices in $u$ not in the same coset, either $\beta_w$ does
not penetrate $w'H$ or $\beta_{w'}$ does not penetrate $wH$, by construction of $\beta$.  Without loss of generality 
we assume the latter.  $\beta_w$ can penetrate $wH$ by at most $c(1)$ so $d(w, w^*) \leq c(1)$.  Then 
$d( w', w'^* ) \leq d( w', w ) + d( w, w^* ) + d( w^*, w'^* ) \leq 1 + c(1) + M$.  The distance from $w$ and $w'$ to
$v$ is bounded by the length of $u$, and the length from $w^*$ and $w'^*$ to $v^*$ is bounded by the length of $u^*$,
so the corresponding wedge between $v$, $v^*$, $w$, $w^*$, $w'$ and $w'^*$ has surface area
also bounded by $5d^1(K')$.
Otherwise, $\beta_v$ does penetrate $gH$ so none of the $\beta_{w_i}$ can penetrate $vH$.  Thus $\beta_v$ can
penetrate $vH$ by at most $2c(1)\left(\lambda(2L+1)+C\right)\left( 2M + 2T + 4 c(1) + 2 \right)$.
Then $d(w_i, w^*_i) \leq d( w_i, v ) + d( v, v^* ) + d( v^*, w^*_i )$ so
these wedges also have surface area bounded by $5 d^1(K'')$. When all of these wedges are included as $3$-cells,
we see that the cones constructed are contractible. 

We adjust the constants so that $K$ is the maximum of what were above called $K$, $K'$, and $K''$. 
The $3$-cells of $\Gamma$ consist of the $3$-cells of $X$ and its translates, as well as a class of mixing
$3$-cells.  The mixing $3$-cells are attached to the $2$-cycles (combinatorial spheres), the 
`drums' of \cite{ger}, consisting of no more than $1 + d^1( (2L+1)J ) + J d^1( 4L+2 ) + 5d^1(K) + W + 2d^1(W)$ cells, 
whose vertices do not all lie in the same translate of $X$. Here $J$ is the maximum number of vertices in a $2$-cell.

\begin{lem}\label{lem2Filling}
$\Gamma$ has Dehn function for filling $2$-boundaries bounded by $f(x) = x^2 \left( d^2_H(x) + x + 1 \right)$.
\end{lem}
\proof
Let $u$ be a finite $2$-subcomplex of $\Gamma$ without boundary, and let $N = |u|$.
Let $u_1$, $u_2$, $\ldots$, $u_m$ denote the connected components of $u$.
If $u_i$ does not pass through any group element vertex, $u_i$ lies entirely within a single translate
of $X$, so $u_i$ can be filled by at most $d^2_H( N )$ $3$-cells.  
Otherwise, $u_i$ passes through at least one group element vertex.  After translating $u$ we may assume that
this vertex is the identity. Let $\sigma'$ be a $2$-cell of $u_i$.  Denote the vertices of $\sigma'$ by
$v'_1$, $v'_2$, $\ldots$, $v'_k$.  If $v'_i$ is a group element vertex, let $v_i = v'_i$, Otherwise, let $v_i$ be
a group element vertex in the same translate at $v'_i$, at a distance no more than $L$ from $v'_i$.
Let $\sigma$ be the cycle $v_1$ to $v_2$ to $\ldots$ to $v_k$ to $v_1$ of length no more than $(2L+1)k$,
through only group element vertices.  The $2$-boundary with top $\sigma$, bottom $\sigma'$ and sides
the cycles $v_i$ to $v'_i$ to $v'_{i+1}$ to $v_{i+1}$ to $v_i$.  This $2$-boundary has area bounded by
$1 + d^1( (2L+1)k ) + k d^1( 4L+2 )$.  If all these vertices are not in the same translate of $X$, 
this corresponds to one $3$-cell.  Otherwise, it can be filled by $\rho = d^2_H( 1 + d^1( (2L+1)J ) + J d^1( 4L+2 ) )$
$3$-cells.  
As $u_i$ is connected, there is $R > 0$ such that the cone from $e$ to $\sigma'$ can be filled 
by at most $RN + Rd^2_H(RN)$ $3$-cells, as constructed above.  $u_i$ can then be filled by
$N\left( RN + Rd^2_H(RN) + \rho \right)$ $3$-cells.  
Repeating this for each $u_i$, we see that $u$ can be filled by $N^2 \left( RN + Rd^2_H(RN) + \rho \right)$ $3$-cells.
\endproof

Suppose that we have constructed the $n$-cells of $\Gamma$, consisting of those $n$-cells lying in translates of $X$ and 
mixing $n$-cells between the translates, with only finitely many orbits of $n$-cells, as above.  

Let $\Delta$ be a mixing $n$-cell. We `comb' $\Delta$ to the basepoint by combing each of the faces and
make the appropriate identifications as collections of $n$-cells.  Bounding the area of each of these `sides', we
uniformly bound the area of each layer of the cone in terms of the $(n-1)$-dimensional Dehn functions of $\Gamma$
and $X$.    The details are as below.

Assume that the vertices of $\Delta$ are $v'_1$, $v'_2$, $\ldots$, $v'_k$, 
$v'_{k+1}=v'_1$.  For each $i$, if $v'_i$ is a group element vertex let $v_i = v'_i$.  Otherwise, there is a group element vertex $v_i$
within an $X$-distance $L$ from $v'_i$.  Connect $v_i$ to $v_j$ via a shortest edge path through group element
vertices whenever $v'_i$ is joined to $v'_j$ by an edge in $\Delta$.  The length of such an edge path is no more than
$\lambda(2L+1)+C$.  Denote the resulting $n$-boundary by $u$.  If the vertices 
$( v'_{i_1}, v'_{i_2}, \ldots, v'_{i_j})$ are a face of $\Delta$, we refer to the corresponding 
$(n-1)$-complex of $u$ as a face of $u$.
For each vertex $v$ of $u$, denote the point where $\beta_v$ enters the coset $vH$ by $v^*$.  $\beta_v$ splits as the concatenation 
$\beta_v = \beta_{v^*} \gamma_v$ with $\gamma_v$ the path traveled by $\beta_v$ from $v^*$ to $v$.
For $v$ and $w$ adjacent vertices of $u$, there is a path connecting $v^*$ and $w^*$ passing through only group element vertices, 
with length no more than $M$.  The distance from $\beta_{v^*}(t)$ to $\beta_{w^*}(t)$ is no more than $M + T + 2c(1)$.
By connecting $\beta_{v^*}(t)$ to $\beta_{w^*}(t)$ by a geodesic path, for each pair of adjacent vertices $v$ and $w$ of $u$,
we obtain successive deformations of $u$, denoted $u_t$, being coned down to the basepoint of $\Gamma$. If $(w_1, w_2, \ldots, w_j)$ are
the vertices of a face of $u$, $(\beta_{w_1^*}(t), \beta_{w_2^*}(t), \ldots, \beta_{w_j^*}(t))$ forms a face of $u_t$ corresponding
to the face $(w_1, w_2, \ldots, w_j)$ of $u$.  For each such $w$, the distance from $\beta_{w^*}(t)$ to $\beta_{w^*}(t+1)$
is no more than $1$.  If $\sigma_t$ and $\sigma_{t+1}$ are corresponding faces of $u_t$ and $u_{t+1}$, joining $\beta_{w^*}(t)$ to
$\beta_{w^*}(t+1)$ for each $w$ we obtain $n$-subcomplexes.  $u_t$, $u_{t+1}$, and these $n$-subcomplexes together results in
an $n$-boundary.  From our earlier remarks there is a $M'_n$ such that each of these $n$-boundaries has no more than $M'$ $n$-cells.  
This finishes our estimates for the bottom portion of the `cone'.

If $v_1$, $v_2$, $\ldots$, $v_k$ are all in the same translate of $X$, say $gX$, then $u$ and $u^*$ (the $n$-complex obtained
by joining $v^*$ to $w^*$ for $v$ and $w$ adjacent vertices of $u$) are both within $gX$.  The $n$-boundary consisting of $u$, $u^*$,
and the $n$-subcomplexes obtained by connecting the corresponding faces of $u$ and $u^*$ within $gX$, can be filled using only the
$(n+1)$-cells from the $gX$ structure.  Otherwise, not all vertices lie in $gX$.  Suppose that $v$ is a vertex of $u$ not lying in
$gX$, and let $( w_1, w_2, \ldots, w_j )$ be a face of $u$.  As above there is a constant $R_n$ such that the $n$-boundary given by the
wedge consisting of $( w_1, w_2, \ldots, w_j )$, $( w_1^*, w_2^*, \ldots, w_j^*)$, the $n$-complexes obtained by connecting their
corresponding faces, as well as connecting $v$ to $v^*$ to each of the corresponding faces (as in the 3 dimensional case above)
consists of at most $R$ $n$-cells.

The $(n+1)$-cells of $\Gamma$ consist of the $(n+1)$-cells in $X$ and its translates as well as the mixing cells,
which correspond to connected $n$-cycles with at most $M'_n + R_n$ $n$-cells.

\begin{cor}
$\Gamma$ has Dehn function for filling $n$-boundaries bounded by $f(x) = x^2 \left( d^n_H(x) + x + 1 \right)$.
\end{cor}
The proof is nearly identical to the proof of Lemma \ref{lem2Filling}, except for the dimension.
These proofs also show the following.
\begin{cor}
$\Gamma$ is acyclic with finitely many $G$-orbits of cells in each dimension.
\end{cor}

This shows that $\Gamma$ is the universal cover of a classifying space of $G$ with finitely many cells in each dimension.
We have thus proven the first part of the following theorem.

\begin{thm}\label{thmNovConj}
Suppose that $G$ is relatively hyperbolic with respect to a collection of finitely many subgroups, each of which is of type $HF^\infty$ and isocohomological, then $G$ is too. Moreover, if each subgroup is further of property RD, then $G$ satisfies the Novikov conjecture.
\end{thm}

The proof for the second part follows from the work of \cite{DS} and \cite{CM} as explained in the introduction.

We remark that the isocohomological property also implies the rational injectivity of the assembly map for the topological $K$-theory of the Lafforgue's algebra $A_{max}(G)$ \cite{laf} by the works of Puschnigg \cite{pu} and Ji-Ogle \cite{Ji-Ogle}. Thus we have

\begin{cor}\label{NovConjAmax}
Suppose that $G$ is relatively hyperbolic with respect to a collection of finitely many subgroups, each of which is of type $HF^\infty$ and isocohomological, then the assembly map for the topological $K$-theory of the Lafforgue's algebra $A_{max}(G)$ is rationally injective. 
\end{cor}

Finally we remark that the class of groups given in the previous theorem is a large class and has a non-trivial intersection with the class of groups that are coarsely embeddable in Hilbert spaces \cite{yu2}. It is not clear
if either class contains the other.

\end{document}